\documentclass{article}
\usepackage{amsthm}
\usepackage{amsfonts}
\usepackage{amssymb}
\usepackage{amsgen}
\usepackage{amsmath}
\usepackage{amsopn}
\usepackage{verbatim}

\theoremstyle{plain}
\newtheorem{thm}{Theorem}[section]
\newtheorem{lem}[thm]{Lemma}
\newtheorem{prop}[thm]{Proposition}
\newtheorem{cor}[thm]{Corollary}
\newtheorem{conj}[thm]{Conjecture}

\theoremstyle{definition}
\newtheorem{prob}[thm]{Problem}
\newtheorem{defn}[thm]{Definition}
\newtheorem{rem}[thm]{Remark}

\newtheorem{eg}[thm]{Example}
\newcommand{\lra}{\longrightarrow}

\newcommand{\wg}{\wedge}
\newcommand{\ub}{\underbrace}
\newcommand{\ra}{\rightarrow}

\newcommand{\zero}{{\bf 0}}
\newcommand{\one}{{\bf 1}}

\newcommand{\bfu}{{\bf u}}
\newcommand{\bfv}{{\bf v}}

\newcommand{\obfj}{{\overrightarrow{\boldsymbol \j}}}

\newcommand{\bfi}{{\bf i}}
\newcommand{\bfj}{{\bf j}}
\newcommand{\bfx}{{\bf x}}
\newcommand{\bfy}{{\bf y}}

\newcommand{\bfk}{{\bf k}}

\newcommand{\caM}{{\mathcal M}}

\newcommand{\F}{{\mathbb F}}
\newcommand{\Z}{{\mathbb Z}}
\newcommand{\R}{{\mathbb R}}
\newcommand{\CP}{{\mathbb CP}}
\newcommand{\Q}{{\mathbb Q}}
\newcommand{\C}{{\mathbb C}}
\newcommand{\V}{{\mathbb V}}

\DeclareMathOperator{\re}{{Re}}
\DeclareMathOperator{\id}{{id}}
\DeclareMathOperator{\im}{{Im}}
\DeclareMathOperator{\GL}{{GL}}
\DeclareMathOperator{\gr}{{gr}}
\DeclareMathOperator{\Aut}{{Aut}}
\newcommand{\gemL}{{\mathfrak L}}
\newcommand{\gemn}{{\mathfrak n}}
\newcommand{\SSS}{{\mathfrak S}}
\newcommand{\om}{{\omega}}
\newcommand{\gam}{{\gamma}}
\newcommand{\gD}{{\Delta}}
\newcommand{\ga}{{\alpha}}
\newcommand{\gl}{{\lambda}}
\newcommand{\gb}{{\beta}}

\newcommand{\gth}{{\theta}}
\newcommand{\gTH}{{\Theta}}
\newcommand{\gt}{{\tau}}
\newcommand{\ol}{\overline}
\newcommand{\rsa}{\rightsquigarrow}
\newcommand{\ot}{\otimes}
\newcommand{\bul}{{\bullet}}
\newcommand{\caF}{{\mathcal F}}
\newcommand{\caV}{{\mathcal V}}
\newcommand{\caO}{{\mathcal O}}

\newcommand{\gO}{{\Omega}}
\newcommand{\gL}{{\Lambda}}
\newcommand{\caS}{{\mathcal S}}
\newcommand{\caL}{{\mathcal L}}
\newcommand{\caD}{{\mathcal D}}
\newcommand{\Lra}{\Longrightarrow}
\newcommand{\lms}{\longmapsto}

\begin{document}
 
\title{Variations of Mixed Hodge Structures of Multiple 
Polylogarithms}
\author{Jianqiang Zhao\footnote{Partially supported by NSF grant DMS0139813}}
\date{}
\maketitle

{\bf Abstract.}
It's well known that multiple polylogarithms give rise to 
good unipotent variations of mixed Hodge-Tate structures.
In this paper we shall {\em explicitly} determine these structures 
related to multiple logarithms and some other multiple polylogarithms 
of lower weights. The purpose of this explicit construction
is to give some important applications: First we study of the limit
mixed Hodge-Tate structures and make a conjecture relating the variations
of mixed Hodge-Tate structures of multiple logarithms to those of
general multiple {\em poly}logarithms. Then following  
Deligne and Beilinson we describe an
approach to defining the single-valued 
real analytic version of the multiple polylogarithms which
generalizes the well-known result of Zagier on
classical polylogarithms. In the process we find some interesting 
identities relating  single-valued multiple polylogarithms of the 
same weight $k$ when $k=2$ and 3. At the end of this paper, 
motivated by Zagier's conjecture we pose
a problem which relates the special values of multiple
Dedekind zeta functions of a number field to the single-valued
version of multiple polylogarithms.


\section{Introduction}
In early 1980s Deligne \cite{D} discovers that the dilogarithm gives 
rise to a good variation of mixed Hodge-Tate structures. 
This has been generalized to polylogarithms (cf.~\cite{Hain}) 
following Ramakrishnan's computation of the monodromy
of the polylogarithms. The monodromy computation also 
yields the single-valued variant $\caL_n(z)$ of the 
polylogarithms (cf.~\cite{BD,Zag}). These functions in turn
have significant applications in arithmetic such as Zagier's conjecture
\cite[p.622]{Zag}. On the other hand,
as pointed out in \cite{G}, ``higher cyclotomy theory'' should 
study the multiple polylogarithm motives at roots of unity, not only
those of the polylogarithms. For this reason we want to look at the 
variations of mixed Hodge structures associated with the 
multiple polylogarithms and see how far we can generalize the 
classical results. In theory such variations of mixed Hodge 
structures are well known to the experts. The purpose of our 
explicit construction is to give some important applications. 

For any positive integer $m_1,\dots,m_n$,
the multiple polylogarithm is defined as follows:
\begin{equation}\label{pexp}
Li_{m_1,\dots,m_n}(x_1,\dots,x_n) =\sum_{0 <
k_1<k_2<\dots<k_n} \frac{x_1^{k_1}x_2^{k_2} \dots
x_n^{k_n}}{k_1^{m_1}k_2^{m_2}\dots k_n^{m_n}},\qquad |x_i|<1.
\end{equation}
We call $n$ the {\em depth} and
$K:=m_1+\dots+m_n$ the {\em weight}. When the depth
$n=1$ the function is nothing but the classical polylogarithm.
More than a century ago H.~Poincar\'e \cite{Po} already knew that 
hyperlogarithms
$$F_n\left({{a_1,\dots, a_n}\atop{b_1,\dots,b_n}}\Big|z\right)
=\int_{b_n}^z\cdots\int_{b_2}^{t_3}\int_{b_1}^{t_2} \frac{d t_1}{t_1-
a_1}
\frac{d t_2}{t_2-a_2}\cdots \frac{d t_n}{t_n-a_n}$$
are important for solving differential equations. 
We observe that although the multiple polylogarithm
can be represented by the iterated path integral in the 
sense of Chen \cite{Ch1}
\begin{equation}\label{old}
Li_{m_1,\dots,m_n}(x_1,\dots,x_n)=(-1)^n
F_K\bigg({{a_1,\overbrace{0,\dots,0}^{m_1-1\text{ times}},
\dots, a_n,\overbrace{0,\dots,0}^{m_n-1\text{ times}}}\atop
{0\ ,\ 0,\dots,0\ ,\dots,\ 0\ \, ,\ 0,\dots,0}}\Big|1\bigg),
\end{equation}
where $a_i=1/(x_i\dots x_n)$ for $1\le i\le n$, it is not
obvious that this actually yields a genuine analytic
continuation in the usual sense when $n\ge 2$. 

According to the theory of framed mixed Hodge-Tate structures
the multiple polylogarithms are period functions
of some variations of mixed Hodge-Tate structures (see \cite{BGSV},
\cite[\S12]{Gicm} and \cite[\S3.5]{Gicm}). Wojtkowiak \cite{Wo} 
studies mixed Hodge structures of iterated integrals over
$\CP^1\setminus\{0,1,\infty\}$ and investigates functional equations
arising from there. In this paper we adopt  
a different approach and compute {\em explicitly} the variations of 
mixed Hodge-Tate structures related to the multiple logarithms
$$\gemL_n(x_1,\dots,x_n):=
Li_{\underbrace{\scriptstyle 1,\dots,1}_{\scriptstyle n\text{ times}}}
(x_1,\dots,x_n).$$
This work relies on our new definition of analytic 
continuation of the multiple polylogarithms given in another paper
\cite{Zana}, by  using Chen's iterated 
path integrals over $\CP^n\setminus D_n$ with some non-normal 
crossing divisor $D_n$. In order to have reasonable variations 
we should be able to control their behavior at ``infinity''. 
This requires us to deal with the natural 
extension of the variations to the infinity using the 
classical result of Deligne \cite[Proposition 5.2]{Del}.
By the same idea we are able to treat all the weight three
multiple polylogarithms and present a result for the double
polylogarithms. From the examples we make the following 
\begin{conj}
The variations of mixed Hodge-Tate structures related to
any multiple polylogarithm can be produced as the variations
of some limit mixed Hodge-Tate structures related to some suitable
choice of multiple logarithm.
\end{conj}
We point out that the old form \eqref{old} of multiple 
polylogarithms is {\em not} suitable for the
investigation of the MHS at the infinity because it is
even not obvious from this form what the ``infinity'' is
exactly.  

As another important application of the our explicit computation, 
in the last section of this paper
we describe an approach to computing the single-valued real analytic
version of the multiple polylogarithms
following an idea of Beilinson and Deligne \cite{BD}. We
find some some interesting identities relating single-valued 
multiple polylogarithms of the same weight $k$ when $k=2$ and 3. For
example, we find the single-valued real analytic
double logarithm (see Eqs. \eqref{llsing} and \eqref{llsing2})
$$\aligned
\caL_{1,1}(x,y)=&\im\big(Li_{1,1}(x,y)\big)
-\arg(1-y)\log|1-x|-\arg(1-xy)\log\Bigl|\frac{x(1-y)}{x-1}\Bigr|\\
=&\caL_2\Bigl(\frac{xy-y}{1-y}\Bigr)
-\caL_2\Bigl(\frac{y}{y-1}\Bigr)-\caL_2(xy)
\endaligned$$
where $\caL_2(z)$ is the famous single-valued dilogarithm.

The motivation of this paper comes from \cite[\S2,3]{G} 
where the Hodge-Tate structures associated with the
double logarithms are discussed, and from \cite{BD} where an
elegant construction of the single-valued real analytic version
of classical polylogarithms are given. The author wishes to thank
R. Hain for answering some of my (perhaps silly) questions concerning 
the good unipotent variations
of mixed Hodge structures. H.~Gangl kindly informed the author of
the preprint \cite{Wo} of Wojtkowiak in which conjectures generalizing
Zagier's are also considered.

As usual HS stands for ``Hodge structure'' and
MHS for ``mixed Hodge structure(s)''.

\section{Multiple logarithms}
We follow the notation in \cite{Zana} in this paper. Recall that
we have an index set
$$\SSS_n=\{\bfi=(i_1,\dots,i_n):\ 0\le i_t\le 1
\text{ for }t=1,\cdots, n\}$$
equipped with a weight function
$$|(i_1,\dots,i_n)|=i_1+\cdots+i_n$$
and two different orderings: a complete ordering $<$
and a partial $\prec$.
If $|\bfi|<|\bfj|$ then $\bfi<\bfj$ (or, equivalently, $\bfj>\bfi$).
If $|\bfi|=|\bfj|$ then the usual lexicographic order from left to right
is in force with $0<1<\cdots $.
The partial ordering is defined as follows.
Let $\bfi=(i_1,\dots,i_n)$
and $\bfj=(j_1,\dots,j_n)$. We set $\bfj\prec \bfi$
(or, equivalently, $\bfi \succ\bfj$) if $j_t\le i_t$ for every 
$1\le t\le n$. For example $(0,0,1,0)\prec (0,1,1,0)$ in 
$\SSS(1,1,1,1)$ but $(1,0,0,0)\not\prec (0,1,1,0)$ and 
$(1,0,0,0)\not\succ (0,1,1,0)$. Clearly
$\bfj\prec \bfi$ implies $\bfj< \bfi$ but not vice versa.

For any $\bfi=(i_1,\dots,i_n)\in \SSS_n$ with $i_s=0$ we define
$$\text{pos}(\bfi,\bfi+\bfu_s)=s$$
as the position where the component is increased by 1.
For example $\text{pos}\big((1,0),(1,1)\big)=2$.
We define the position functions $f_n^1,\dots,f_n^n$ on
$\obfj\in \SSS_n^n$ as follows:
$$f_n^1(\obfj)=1,\quad
f_n^t(\obfj)=\text{pos}\big(\bfj_{t-1},\bfj_t\big),
\text{ for \ } 2\le t\le n.$$ 
These functions tell us the places where the
increments occur in the queue of $\obfj$. Set
$$w_1(\bfx):=d\log\Bigl(\frac{1}{1-x_1}\Bigr);\quad
w_t(\bfx):=d\log\Bigl(\frac{1-x_{t-1}^{-1}}{1-x_t}\Bigr),
\text{ for \ }2\le t\le n.$$
The analytic 
\begin{prop} {\em (\cite[Proposition 5.1]{Zana})}
The multiple logarithm $\gemL_n(\bfx)$ is a
multi-valued holomorphic function on 
$$\prod_{1\le j\le n}(1-x_j)
\prod_{1\le j<k\le n} \big(1-x_j\dots x_k\big)=0.$$
and can be expressed by
\begin{equation} \label{formula}
\gemL_n(\bfx)=\sum_{\obfj=(\bfj_1,\dots,\bfj_n)\in\SSS_n^n}
\int_{\zero}^{\bfx} w_{f_n^1(\obfj)}(\bfx(\bfj_1))
w_{f_n^2(\obfj)}(\bfx(\bfj_2))
\cdots w_{f_n^n(\obfj)}(\bfx(\bfj_n)),
\end{equation}
where the path from $\zero$ to $\bfx$ lies in $S'_n$.
\end{prop}

\section{Multiple logarithm variations of MHS}\label{genvar}
In this section we will define the variation matrix $\caM_{[n]}(\bfx)$
coming from the multiple logarithms of depths up to $n$.
We will show that it is a $2^n\times 2^n$ multi-valued matrix which
defines a good variation of a MHS over $S_n=\C^n\setminus D_n$
where $D_n$ is the divisor defined by
$$\prod_{1\le j\le n}x_j(1-x_j)
\prod_{1\le j<k\le n} \big(1-x_j\dots x_k\big)=0.$$

\begin{rem} In fact, the irreducible component $x_n=0$ in $D_n$ is
not needed in the case of multiple logarithms. But the variation
matrix corresponding to general
multiple polylogarithms may have singularities along this component,
for example, $\caM_{1,2}(x_1,x_2)$ of the double polylogarithm 
$Li_{1,2}(x_1,x_2)$. See chapter~\ref{wt3var}.
\end{rem}

\subsection{Definition of variations of MHS: a review}
In this section we briefly review  the theory of
variations of MHS.

A pure ($\Q$-)HS of weight $k$ consists of a finitely generated
abelian group $H(\Z)$ and a decreasing Hodge filtration $\caF^\bul$ on
$H(\C):=H(\Z)\ot_\Z \C$ such that $H(\C)=\caF^p\oplus \ol{\caF^{k-p+1}}$
for all integers $p$. Here the ``bar'' is the complex conjugation on the
second factor of the tensor product. A special example is the Tate
structure $\Z(-k)$ of weight $2k$ consists of $H(\Z)=\Z$ and the
filtration $\caF^p=0$ for $p>k$ and $\caF^p =H(\C)$ for $p\le k$.

A MHS consists of a finitely generated
abelian group $H(\Z)$ and two filtrations: an increasing
weight filtration $W_\bul$ on $H(\Q):=H(\Z)\ot_\Z \Q$ and
a decreasing filtration $\caF^\bul$ on $H(\C)$,
which are compatible in the following sense.
On each graded piece of the weight filtration $\gr_k^W=W_k/W_{k-1}$
the induced Hodge filtration determined by
$$\caF^p(\gr_k^W)(\C)=\text{\raisebox{1ex}{$(\caF^p\cap W_k(\C) + 
W_{k-1}(\C) )$}}\big/
\text{\raisebox{-1ex}{$W_{k-1}(\C)$}} $$
is a pure Hodge structure of weight $k$ where $W_k(\C):=W_k\ot_Z \C$.
If all the pure Hodge structures induced as above are direct sums
of Tate structures then we say the MHS is a Tate structure.
For a mixed Hodge-Tate structure we can put a framing as in
\cite[\S1.3.4, \S1.4]{BGSV}.

Following Steenbrink and Zucker \cite[Definitions 3.1, 3.2 and 3.4]{SZ} 
we have
\begin{defn}
A {\em variation of HS of weight $k$} defined over
$\Q$ and a complex manifold $S$ is
a collection of data $(\V_\Q,\caF^\bul)$ where

(a) $\V_\Q$ is a locally constant sheaf (local system)
of $\Q$-vector spaces on $S$,

(b) $\caF^\bul$ is a decreasing filtration by holomorphic subbundles of 
the
locally free sheaf $\caV=\caO_S\ot_\Q \V_\Q$.

(c) At each $s\in S$, $\caF^\bul$ induces the Hodge filtration 
$\caF^\bul_s$
of a Hodge structure of weight $k$ on the fiber $\caV_s$ of $\caV$ such
that
\begin{quote}
\noindent
(i) whenever $p+q=k$ one has $\caV_s=\caF^p_s\oplus \ol{\caF^{q+1}_s},$
where the ``bar'' denotes the complex conjugation,

\noindent
(ii) equivalently, one has $\caV_s=\bigoplus_{p+q=k} H_s^{p,q}$ where  
$H_s^{p,q}=\caF^p_s\cap \ol{\caF^q_s}.$
\end{quote}

(d) (Griffiths transversality)
Under the connection $\nabla$ in $\caV$,
$$\nabla \caF^p\subset \gO_S^1\ot_{\caO_S} \caF^{p-1}\qquad\text{ for 
all } 
p.$$
\end{defn}

\begin{defn} {\em A polarization over} $\Q$ of a  variation of Hodge
structure of weight $k$ over $\Q$ is a non-degenerated and flat
bilinear pairing:
$$\gb:\V_\Q\times \V_\Q\lra \Q,$$
such that $\gb$ is $(-1)^k$-symmetric, and the Hermitian form 
$\gb_s(C_sv, \bar w)$ is positive on each fiber. Here $C_s$ denotes
the Weil operator with respect to $\caF_s$, namely the direct sum of
multiplication by $i^{p-1}$ on $H_s^{p,q}$. A variation is called
{\em polarizable} (over $\Q$) if it admits a polarization (over $\Q$).
\end{defn}

\begin{defn}\label{varHSdef}
A {\em variation of MHS} defined over
$\Q$ and a complex manifold $S$ is
a collection of data $(\V_\Q,W_\bul,\caF^\bul)$ where

(a) $\V_\Q$ is a local system of $\Q$-vector spaces on $S$,

(b) $W_\bul$ is an increasing filtration of the $\V_\Q$ by local 
subsystems,

(c) $\caF^\bul$ is a decreasing filtration by holomorphic subbundles of
$\caV=\caO_S\ot_\Q \V_\Q$.

(d) $\nabla \caF^p\subset \gO_S^1\ot_{\caO_S} \caF^{p-1}$ for all $p$.

(e) The data 
\begin{equation}\label{fibre}
\big(\gr_k^W \V_\Q, \caF^\bul\big(
\text{\raisebox{.5ex}{$\caO_S\ot_E W_k$}}\big/
\text{\raisebox{-.5ex}{$\caO_S\ot_E W_{k-1}$}}\big)\big)
\end{equation}
is a variation of HS of weight $k$ defined over $\Q$; or equivalently,
on the fiber over $s\in S$, $(V_s,W_s,\caF_s)$ is a MHS defined over 
$\Q$.

(f) If the induced collection of variations of HS \eqref{fibre} are all
polarizable then the MHS is called graded-polarizable.
\end{defn}

\begin{rem} By extension of scalars in $\V_\Q$ one can 
define $\V_\F$ for any field $\F$ such that $\Q\subset \F\subset \R$. 
\end{rem}

Giving a local system $\V_\Q$ is
equivalent to specifying its monodromy representation
$$\rho_x: \pi_1(S,x)\lra \Aut_\Q \caV_x.$$
A variation is called {\em unipotent} if this representation
is unipotent. From Proposition 1.3 of \cite{HZ2} we know that
a variation of MHS $(\V_\Q,W_\bul,\caF^\bul)$
is unipotent if and only if each of the variations of Hodge
structure $\gr_k^W \V_\Q$ is constant.

In general, the behavior of a variation of MHS over a non-compact
base $S$ at ``infinity'' is very hard to control. Steenbrink 
and Zucker \cite{SZ}
consider the case when $S$ is a curve and define the admissibility
condition at infinity. For higher dimensional $S$, Kashiwara, M. Saito, 
and others define a variation
over $S$ to be admissible if its restriction to
every curve is admissible in the sense of Steenbrink-Zucker.

However, the behavior of {\em unipotent} variations of MHS
at infinity can be controlled rather easily. We have the 
classical result of Deligne \cite[Proposition 5.2]{Del}
which defines the {\em canonical extension} $\tilde\caV$ of $\caV$.

\begin{thm}\label{del} {\em (Deligne)} 
Let $\tilde{S}$ be a normalization of $S$. Let
$(\V_\Q,W_\bul,\caF^\bul)$ be a unipotent variation of MHS over $S$
whose associated connection in $\caV$ is integrable. Then 

\noindent(a) There
is a unique extension $\tilde\caV$ of $\caV$ over $\tilde{S}$ 
satisfying the following equivalent conditions:
\begin{quote}
(i) Inside every section of $\tilde\caV$, every flat 
section of $\caV$ increases at most at the rate of 
$O(\log^k|\!|x|\!|)$ ($k$ large enough)
on every compact set of $D=\tilde{S}-S$.

(ii) Similarly, every flat section of $\caV^{\vee}$ (the dual)
increases at most at the rate of $O(\log^k|\!|x|\!|)$ 
($k$ large enough).
\end{quote}

\noindent(b) The combination of the two conditions (i) and (ii) is
equivalent to the combination of the following two conditions:
\begin{quote}
(iii) In terms of any local basis of $\tilde\caV$ the connection matrix
$\boldsymbol \om$ of $\caV$ has at most logarithmic singularities along 
$D$.

(iv) The residue of $\boldsymbol \om$ along any irreducible component 
of $D$ is
nilpotent.
\end{quote}
\end{thm}
We will verify conditions (iii) and (iv) by Proposition \ref{griff}
for the multiple logarithm variations of MHS. They are unipotent
variations by Theorem~\ref{thm:mono}.

\begin{defn}
Let $\tilde{S}$ be a compactification of $S$.
Then a unipotent variation of MHS
$(\V_\Q,W_\bul,\caF^\bul)$ over $S$ is said to be {\em good} if it
satisfies the following conditions at infinity
\begin{quote}
(1) the Hodge filtration bundles $\caF^\bul$
extend over $\tilde{S}$ to sub-bundles $\tilde{\caF}^\bul$ of
the canonical extension $\tilde{\caV}$ of $\caV$ such that they induce
the corresponding thing for each pure subquotient $\gr_k^W \V_\Q$,

(2) for the nilpotent logarithm $N_j$ of a
local monodromy transformation about a component $D_j$ of $D$, the
weight filtration of $N_j$ relative to $W_\bul$ exists.
\end{quote}
\end{defn}

A slightly different definition first appeared in \cite{HZ1,HZ2}
with the extra assumption that $D=\tilde{S}-S$ is a normal crossing 
divisor.
In these papers Hain and Zucker classified good unipotent
variations of MHS on algebraic manifolds. With constant pure weight 
subquotients these variations behave well at infinity. 

\subsection{The variation matrix} 
The double logarithm was treated in
\cite[\S2]{G} by Goncharov. We noticed an apparent typo that the term
$2\pi i\log x$ in the matrix $A_{1,1}(x,y)$ on page 620 should be 
replaced
by $2\pi i\log (1-x)$.
We first rewrite $A_{1,1}(x,y)$ as
$\caM_{1,1}(x,y)$ below because we will use induction
starting from this form of double logarithm variation of MHS in several
proofs later.
\begin{equation}\label{n=2}
\caM_{1,1}(x,y)=\left[\begin{matrix}
1&0&0&0\\
\gemL_1(y)  & 2\pi i&0&0\\
\gemL_1(xy) & 0&2\pi i&0\\
\gemL_2(x,y)& 2\pi i\gemL_1(x)&
	2\pi i\gemL_{1}\big(\frac{1-xy}{1-x}\big)&(2\pi i)^2
\end{matrix}\right].
\end{equation}
This is essentially the same as defined in \cite{G} up to signs.

We now begin to define the variation matrix $\caM_{[n]}(\bfx)$ for
every $\bfx\in S_n$.  
\begin{defn}\label{def:Eij}
For $1\le s\le n$ write as before
$$\gth_s=\gth_s(\bfx)=\frac{dt}{t-a_s(\bfx)}=\frac{dt}{t-a_s}.$$
Suppose $|\bfi|=k$ and $i_{\gt_1}=\dots=i_{\gt_k}=1$.
Suppose $|\bfj|=l$ and $j_{t_1}=\dots=j_{t_l}=1$.
\begin{description}
\item[{\sl{\em (1)}}] If $\bfj\not\prec\bfi$, we define the
$(\bfi,\bfj)$-th entry of $\caM_{[n]}(\bfx)$ to be $0$.
\item[{\sl{\em (2)}}] If $\bfj\prec \bfi$ then we let
$t_r=\gt_{\ga_r}$ for $1\le r\le l$ and 
\begin{equation}\label{yss}
\bfx(\bfi)=\bfy=(y_1,\dots,y_k),\qquad
y_m=\prod_{\ga=\gt_m}^{\gt_{m+1}-1}x_\ga
	=\frac{a_{\gt_{m+1}}(\bfx)}{a_{\gt_m}(\bfx)}, 
	\quad 1\le m\le k.
\end{equation}
with $\gt_{k+1}=n+1$ and $a_{n+1}=1$. Set $t_0=\ga_0=0$,
$t_{l+1}=n+1$, $\ga_{l+1}=k+1$, $a_0(\bfx)=a_0(\bfy)=0$
and $a_{n+1}(\bfx)=a_{k+1}(\bfy)=1$.
Define the $(\bfi,\bfj)$-th entry of $\caM_{[n]}(\bfx)$
as $(2\pi i)^lE_{\bfi,\bfj}(\bfx)$ where
\begin{equation}\label{coldefn}
\aligned
E_{\bfi,\bfj}(\bfx)=\gam_{\rho_\bfi(\bfj)}^k(\bfy)
:=&(-1)^{k-l}\prod_{r=0}^l
	\int_{a_{\ga_r}(\bfy)}^{a_{\ga_{r+1}}(\bfy)}
	\gth_{\ga_r+1}(\bfy)\cdots\gth_{\ga_{r+1}-1}(\bfy)\\
=&(-1)^{k-l}\prod_{r=0}^l \int_{p_r}
\gth_{\gt_{\ga_r+1}}(\bfx)\cdots \gth_{\gt_{\ga_{r+1}-1}}(\bfx)
\endaligned
\end{equation}
Here the $l+1$ paths $p_0,\dots,p_l$ for the $l+1$
integrals are {\em independent of } $\bfi$ where $p_r$ is
any fixed contractible path from $a_{t_r}$ to $a_{t_{r+1}}$ in the 
punctured
complex plane $\C\setminus \bigcup_{t_r<s<t_{r+1}} \{a_s\}$,
and the integral $\int_{p_r}=1$ if $\ga_r+1=\ga_{r+1}.$ 
We get the second equality by observing that
$$a_m(\bfy)=(y_m\dots y_k)^{-1}=a_{\gt_m}(\bfx)\Lra
a_{\ga_r}(\bfy)=a_{\gt_{\ga_r}}(\bfx)=a_{t_r}(\bfx).$$
\end{description}
\end{defn}

\begin{prop}\label{intexp} Suppose $\bfi$ and $\bfj$ are given
as in Definition~{\em \ref{def:Eij}(2)}. As multi-valued functions
\begin{align}\label{Eij}
E_{\bfi,\bfj}(\bfx)=&\prod_{r=0}^l\gemL_{\ga_{r+1}-\ga_r-1}\left(
\frac{a_{\gt_{\ga_r+2}}(\bfx)-a_{t_r}(\bfx)}
{a_{\gt_{\ga_r+1}}(\bfx)-a_{t_r}(\bfx)}, \cdots,
\frac{a_{t_{r+1}}(\bfx)-a_{t_r}(\bfx)}
{a_{\gt_{\ga_{r+1}-1}}(\bfx)-a_{t_r}(\bfx)}\right)\\
=&\gemL_{\ga_1-1}\big(x_{\gt_1}\cdots x_{\gt_2-1},
	x_{\gt_2}\cdots x_{\gt_3-1},\cdots,
	 x_{\gt_{\ga_1-1}}\cdots x_{t_1-1}\big)\cdot \notag\\
\ &\cdot\prod_{r=1}^l\gemL_{\ga_{r+1}-\ga_r-1}\left(
\frac{1-x_{t_r}\cdots x_{\gt_{\ga_r+2}-1}}
	{1-x_{t_r}\cdots x_{\gt_{\ga_r+1}-1}}, \cdots,
\frac{1-x_{t_r}\cdots x_{t_{r+1}-1}}
	{1-x_{t_r}\cdots x_{\gt_{\ga_{r+1}-1}-1}}\right).\label{Eijlog}	
\end{align}
Here $\gemL_0=1$ and $a_0=0$.
\end{prop}
\begin{proof}
By direct and simple calculation we get
\begin{multline*}
(-1)^{\ga_{r+1}-\ga_r-1}\int_{p_r}
\gth_{\gt_{\ga_r+1}}(\bfx)\cdots \gth_{\gt_{\ga_{r+1}-1}}(\bfx)\\
=\gemL_{\ga_{r+1}-\ga_r-1}\left(
\frac{a_{\gt_{\ga_r+2}}(\bfx)-a_{t_r}(\bfx)}
{a_{\gt_{\ga_r+1}}(\bfx)-a_{t_r}(\bfx)},
\frac{a_{\gt_{\ga_r+3}}(\bfx)-a_{t_r}(\bfx)}
{a_{\gt_{\ga_r+2}}(\bfx)-a_{t_r}(\bfx)}, \dots,
\frac{a_{t_{r+1}}(\bfx)-a_{t_r}(\bfx)}
{a_{\gt_{\ga_{r+1}-1}}(\bfx)-a_{t_r}(\bfx)}\right).
\end{multline*}
The proposition follows immediately.
\end{proof}

\begin{eg} On the last row of $\caM_{[n]}(\bfx)$ one has
\begin{align}\label{gamj}
E_{\one,\bfj}(\bfx)=&\gam_{\bfj}^n(\bfx)
	=\prod_{r=0}^l \gemL_{t_{r+1} - 
	t_r -  1}
	\Bigl(\frac{a_{t_r+2} -  a_{t_r}}
	{a_{t_r+1} -  a_{t_r}},\cdots,
	\frac{a_{t_{r+1}} -  a_{t_r}}
	{a_{t_{r+1}-1} -  a_{t_r}}\Bigr)\\
	=&\prod_{r=0}^l \gemL_{t_{r+1} - 
	t_r -  1}
	\Bigl(\frac{1-x_{t_r}x_{t_r+1}}{1-x_{t_r}},\cdots,
	\frac{1-x_{t_r}\cdots x_{t_{r+1}-1}}
	{1-x_{t_r}\cdots x_{t_{r+1}-2}}\Bigr)\notag
\end{align}
where $\gemL_0=1$ and $x_0=\infty$. In particular,
$E_{\one,\zero}=\gam_{\zero}^n(\bfx)=\gemL_n(\bfx)$
and $E_{\one,\one}=\gam_{\one}^n(\bfx)=1$.
\end{eg}

We now fix a standard basis $\{e_\bfi: \bfi\in \SSS_n\}$
of $\C^{2^n}$ consisting of column vectors.
Suppose $|\bfi|=k$. It follows from definition that
the $\bfi$-th row is
\begin{equation}\label{kthblock}
R_\bfi:=\sum_{\bfj\prec\bfi} (2\pi i)^{|\bfj|}
	\gam_{\rho_\bfi(\bfj)}^k\big(\bfx(\bfi)\big)e_\bfj^T
=(2\pi i)^k e_\bfi^T+\sum_{\bfj\precneqq \bfi} (2\pi i)^{|\bfj|}
	\gam_{\rho_\bfi(\bfj)}^k\big(\bfx(\bfi)\big) e_\bfj^T
\end{equation}
where $e_\bfj^T$ are now row vectors.
Note that $\gam_{\rho_\bfi(\bfi)}^k=\gam_{\one_k}^k=1$ by definition.
It is clear that the first entry (i.e. $\bfj=\zero$)
of this row is $\gemL_k \big(\bfx(\bfi)\big)$.

Let us call the minor of $\caM_{[n]}(\bfx)$ consisting of
rows beginning with $k$-tuple logarithms 
{\em the $k$-th block}. It has ${n \choose k}$ rows with row
indices $|\bfi|=k$.

\begin{lem}\label{kthdiag}
The matrix $\caM_{[n]}(\bfx)$ is a lower triangular matrix. Moreover,
the columns with $|\bfj|=k$ of the $k$-th block of $\caM_{[n]}(\bfx)$
is $(2\pi i)^k$ times the identity matrix of rank ${n \choose k}$.
\end{lem}
\begin{proof} The lemma follows directly from equation \eqref{kthblock}
because if $\bfj\precneqq \bfi$ then $\bfj<\bfi$.
\end{proof}

\begin{lem}\label{jthcolumn} The $\bfj$-th column of $\caM_{[n]}(\bfx)$ 
is
$$(2\pi i)^{|\bfj|}C_\bfj=(2\pi i)^{|\bfj|}\sum_{\bfi \succ\bfj}
  \gam_{\rho_\bfi(\bfj)}^{|\bfi|}\big(\bfx(\bfi)\big) e_\bfi$$
where $\bfx(\bfi)$ are defined by equation \eqref{yss} depending on 
$\bfi$.
\end{lem}
\begin{proof} Use equation \eqref{kthblock}.
\end{proof}

\begin{eg} By definition or the above proposition the first column
$$C_\zero(\bfx)=\big[\gemL_{|\bfi|}(\bfx(\bfi)):\ \bfi\in \SSS_n\big]^T$$
where $\gemL_0=1$. 
\end{eg}

\begin{prop} \label{griff}
The columns of $\caM_{[n]}(\bfx)$ form the set of the
fundamental solutions of the following
system of differential equations
\begin{equation} \label{genDE}
\left\{\aligned
d\,X_\zero=&0,\\
d\,X_\bfi=&\sum_{|\bfk|=|\bfi|-1,\, \bfk\prec \bfi} X_\bfk \,
d\,\gam_{\rho_\bfi(\bfk)}^{|\bfi|}\big(\bfx(\bfi)\big)\quad
	\text{for all } 1\le |\bfi|\le n
\endaligned\right.
\end{equation}
where $\bfx(\bfi)$ is determined as in equation \eqref{yss}.
\end{prop}
\begin{proof} We prove the proposition by induction on $n$.
It is easy to see the proposition is valid for $n=1$ and $n=2$.
We assume that $n\ge 3$ and the proposition is true for $\le n-1$.
Let us now look at the $\bfj$-th column as expressed
in Lemma~\ref{jthcolumn}. The cases $|\bfi|=1$ or
$\bfj>\bfi$ are obvious. Suppose

(1) $1<|\bfi|<n$ and $\bfj\le \bfi$.
There are two cases. (i) $\bfj\not\prec \bfi$. This is trivial
because each term of both sides is zero.
(ii) $\bfj\prec \bfi$. Then there is a $t$ such that $i_t=j_t=0$.
We denote $\bfi'\in \SSS_{n-1}$ the corresponding index
after deleting the $i_t$-th component. By induction
$$\sum_{|\bfk'|=|\bfi'|-1,\, \bfj'\prec \bfk'\prec \bfi'}
\gam_{\rho_{\bfk'}(\bfj')}^{|\bfk'|}\big(\bfx'(\bfk')\big) \,
d\,\gam_{\rho_{\bfi'}(\bfk')}^{|\bfi'|}\big(\bfx'(\bfi')\big)=
d\,\gam_{\rho_{\bfi'}(\bfj')}^{|\bfi'|}\big(\bfx'(\bfi')\big)
$$
where we set $\bfx'=(x_1,\dots,x_{i_t-1},x_{i_t}x_{i_t+1},
x_{i_t+2},\dots,x_n)$.
Since $|\bfi'|=|\bfi|$ and $|\bfk'|=|\bfk|$ we can get
the desired equation by inserting $0$ before the $i_t$-th
components of $\bfi'$, $\bfj'$ and $\bfk'$, i.e., using
the embedding $\iota_{i_t}$.

(2) $\bfi=\one$ and $|\bfj|=l$. We need to show
\begin{equation}\label{lastrow}
d\gam_\bfj^n(\bfx)=\sum_{|\bfk|=n-1,\, \bfj\prec \bfk}
\gam_{\rho_\bfk(\bfj)}^{n-1}\big(\bfx(\bfk)\big)\,
d\gam_\bfk^n(\bfx).
\end{equation}
This is trivial when $l=n$. The case $l=0$ follows from 
$$ d\gemL_n(\bfx)=\sum_{t=1}^n
\gemL_{n-1}(x_1,\dots,x_{t-2},x_{t-1}x_t,x_{t+1},\dots,x_n)\,
d\log\frac{1-x_{t-1}^{-1}}{1-x_t}.$$
So we may assume $0<l<n$,
$j_{t_1}=\dots=j_{t_l}=1$ and $j_t=0$ for all other indices $t$.
By definition \eqref{gamj} we have
$$\gam_\bfj^n(\bfx)= \sum_{r=0}^l\sum_{t_r<s<t_{r+1}}
\gam_{\rho_{\bfv_s}(\bfj)}^{n-1}\big(\bfx(\bfv_s)\big)\,
d\gam_{\bfv_s}^n(\bfx)$$
where $t_0=0$, $t_{l+1}=n+1$ and
$${\setlength\arraycolsep{1pt}
\begin{array}{rl}
\bfv_s=&(1,\dots,1,\dots\dots\dots,1,0,1,\dots\dots\dots,1,\dots,1).\\
\     &\phantom{(1,\dots,}\uparrow \phantom{,\dots\dots\dots,1}\uparrow \phantom{1,\dots\dots\dots,}\uparrow \\
\ &\phantom{(1,}t_r\text{-th place}  \phantom{,\dots,}s\text{-th place}  
	\phantom{\dots,} t_{r+1}\text{-th place}
\end{array}}$$
Under the retraction map $\rho_{\bfv_s}$ the numbering of the
indices changes as follows: $t\rsa t$ if $t<s$ and
$t\rsa t-1$ if $t>s$. We also have
$$a_t\big(\bfx(\bfv_s)\big)=\begin{cases}
a_t(\bfx) \qquad &\text{if }t<s,\\
a_{t+1}(\bfx)  &\text{if }t>s.\end{cases}$$
Hence for each $s$ such that $t_r<s<t_{r+1}$ the integral expression of
$\gam_{\rho_{\bfk}(\bfj)}^{n-1}\big(\bfx(\bfk)\big)$
is unchanged under $\rho_\bfk$ ($\bfj\prec \bfk$)
except the $\bfv_s$-term. Equation~\eqref{lastrow} now follows
immediately from Leibniz rule and so the proposition is proved.
\end{proof}

\subsection{Monodromy of $\caM_{[n]}(\bfx)$} 
Fix an embedding $\C^n\hookrightarrow \CP^n$. Let $\caD_n=D_n\cup 
(\CP^n\setminus \C^n)$. Let $M_r(\C)$ be the set of $r\times r$
matrices over $\C$.  Put
\begin{equation}
{\boldsymbol \om}=\Bigl( c_{\bfi,\bfj}\Bigr)_{\bfi,\bfj\in \SSS_n}\in
H^0(\CP^n,\Omega^1_{\CP^n}(\log(\caD_n)))\otimes M_{2^n}(\C)
\end{equation}
where
$$c_{\bfi,\bfj}=\begin{cases}
d\gam_{\rho_\bfi(\bfj)}^{|\bfi|}\big(\bfx(\bfi)\big) \quad &\text{if }
|\bfj|=|\bfi|-1,\, \bfj\prec \bfi, \\
0 &\text{otherwise}.
\end{cases}$$
All of the 1-forms in $\boldsymbol \om$
have logarithmic singularity on $\caD_n$ because of
the following. Let $|\bfi|=l$ and $i_{t_1}=\dots=i_{t_l}=1$.
Let $j_{t_s}=0$ so that $|\bfj|=l-1$ and $\bfj\prec \bfi$. Let
$\bfx(\bfi)=\bfy=(y_1,\dots,y_l)$. By definition \eqref{coldefn}
$$\aligned
\gam_{\rho_\bfi(\bfj)}^{|\bfi|}\big(\bfx(\bfi)\big)
\ &=-\int_{a_{s-1}(\bfy)}^{a_{s+1}(\bfy)}
	\gth_s(\bfy)
=-\log\left(\frac{a_{s+1}(\bfy)-a_s(\bfy)}
{a_{s-1}(\bfy)-a_s(\bfy)}\right)\\
\ &
\begin{array}{ll}
=\begin{cases}
-\log(1-y_1)   \\
-\log\Bigl( \frac{\textstyle y_{s-1}(y_s-1)}{\textstyle 1-y_{s-1}}
\Bigr) \quad \end{cases}
\ &\aligned \text{if } &s=1\\
	\text{if } &s\ge 2 \endaligned \\
=\begin{cases}
-\log(1-x_1\dots x_{t_1}) \\
-\log\Bigl(
\frac{\textstyle x_{t_{s-1}}\dots x_{t_s-1}(x_{t_s}\dots x_{t_{s+1}-1}-
1)}
{\textstyle 1-x_{t_{s-1}}\dots x_{t_s-1}}\Bigr)
\end{cases}
\ &\aligned \text{if } &s=1\\
	\text{if } &s\ge 2. \endaligned 
\end{array}
\endaligned$$
\begin{eg} When $n=2$ we have
$${\boldsymbol \om}=\left[\begin{matrix}
0 & \\
-d\log(1-y) &0\\
-d\log(1-xy) & \ &0\\
0 & -d\log(1-x) &-d\log\frac{x(1-y)}{x-1} & 0
\end{matrix}\right]$$
\end{eg}

We proved in Proposition \ref{griff} that $\caM_{[n]}(\bfx)$
is a fundamental solution of first order linear partial differential 
equation
\begin{equation}\label{pde}
d\gL = {\boldsymbol \om}\gL
\end{equation}
where $\gL$ is a possibly multi-valued function $S\lra M_{2^n}(\C)$. 
Moreover $\caM_{[n]}(\bfx)$ is a unipotent matrix for very $\bfx\in S$.
Applying $d$ on equation \eqref{pde} and plugging in 
$\gL=\caM_{[n]}(\bfx)$ we get
$$0=d {\boldsymbol \om}\caM_{[n]}(\bfx)
-{\boldsymbol \om}\wg d\caM_{[n]}(\bfx)
=(d {\boldsymbol \om}-{\boldsymbol \om}\wg{\boldsymbol \om}).
\caM_{[n]}(\bfx)$$
Because $\caM_{[n]}(\bfx)$ is invertible and 
$\boldsymbol \om$ is closed we get
\begin{equation}\label{integrabl}
d {\boldsymbol \om}=0, \quad {\boldsymbol \om}\wg {\boldsymbol \om}=0.
\end{equation}
This shows that ${\boldsymbol \om}$ is integrable.

The main goal of this chapter is to show that if we analytically 
continue 
every integral entry of $\caM_{[n]}(\bfx)$ along a same
loop $q\in \pi_1(S_n,\bfx)$, the resulting
matrix will still be a fundamental solution 
$\caM_{[n]}(\bfx)M(q)$ of \eqref{pde} where $M(q)\in \GL_{2^n}(\Z)$.
In the following we also denote this action of $q$ by $\gTH(q)$ 
operating
on the left. We then define the monodromy representation
$$\aligned
\rho_\bfx: \pi_1(S_n,\bfx)&\lra \GL_{2^n}(\Z)\\
	q&\lms M(q)^T.
\endaligned$$
Here we take the transpose to ensure $\rho_\bfx$ to be a homomorphism
because $M(pq)=M(q)M(p)$ by our convention. From the explicit 
computation in Theorem \ref{thm:mono}
we will see that $\rho_\bfx$ is a unipotent representation.

\begin{thm} \label{thm:mono}
Let $\caM_{[n]}(\bfx)=\big[
E_{\bfi,\bfj}(\bfx)\big]_{\bfi,\bfj\in \SSS_n}$ where
$E_{\bfi,\bfj}(\bfx)$ are defined by
Proposition~{\em \ref{intexp}}. Let $1\le i\le j\le n$ and
$q_{ij}\in \pi_1(S_n,\bfx)$ (resp. $1\le j<n$ and
$q_{j0}$) enclose $\caD_{ij}=\{x_i\dots x_j=1\}$, (resp.
$\caD_{j0}=\{x_j=0\}$) 
only once but no other irreducible component of $D_n$ such that
$\int_{q_{ij}} d\log(1-x_i\dots x_j)=2\pi i$ 
(resp. $\int_{q_{j0}} d\log x_j=2\pi i$). Then
$$
M(q_{j0})=I+\big[n_{\bfi,\bfj}\big]_{\bfi,\bfj\in \SSS_n},\quad
M(q_{ij})=I+\big[m_{\bfi,\bfj}\big]_{\bfi,\bfj\in \SSS_n}$$
where $I$ is the identity matrix of rank $2^n$,
\begin{equation}\label{nij}
n_{\bfi,\bfj}=\begin{cases} 
-1 \quad &\text{if }t_r\le j\le t_{r+1}-2,\ r\ge 1,\
	\bfi=\bfj+\bfu_{s+1}\text{ and }j\le s\le t_{r+1}-2\\
0 \quad &\text{otherwise,}\end{cases}
\end{equation}
and
\begin{equation}\label{mij}
m_{\bfi,\bfj}=\begin{cases}
1 &\text{if $t_r=i\le j\le t_{r+1}-2,\ r\ge 1$,
	$\bfi=\bfj+\bfu_{j+1}$}\\
-1 &\text{if $t_r+1\le i\le j=t_{r+1}-1,\ r\ge 0$,
	$\bfi=\bfj+\bfu_i$}\\
0 \quad &\text{otherwise.}\end{cases}
\end{equation}
Here $\bfi$ and $\bfj$ in the case of $m_{\bfi,\bfj}=\pm 1$ 
and $n_{\bfi,\bfj}=-1$
satisfy the condition in Definition~{\em \ref{def:Eij}(2)}.
\end{thm}
\begin{proof}. 
By definition it is clear that if $\bfi\not\succ \bfj$ then
$\gTH(q) E_{\bfi,\bfj}(\bfx)=E_{\bfi,\bfj}(\bfx)$ which is either 0 or 
1.
Thus we are only concerned with $E_{\bfi,\bfj}$ with $\bfi\succ \bfj $.

We now fix some $\bfj$. If $|\bfj|=n$ then clearly 
$(\gTH(q)-I)C_\one=[0,\dots,0]^T$ for any loop $q$.
This proves the proposition for $|\bfj|=n$. We now assume $|\bfj|<n$. 
Let $\bfi$ and $\bfj$ be given as in Definition \ref{def:Eij}(2).
By equation \eqref{Eij} 
\begin{multline*}
E_{\bfi,\bfj}(\bfx)=\gemL_{\ga_1-1}\big(x_{\gt_1}\cdots 
	x_{\gt_2-1},x_{\gt_2}\cdots x_{\gt_3-1},\cdots,
	 x_{\gt_{\ga_1-1}}\cdots x_{t_1-1}\big)\cdot \\
\cdot\prod_{r=1}^l\gemL_{\ga_{r+1}-\ga_r-1}\left(
\frac{1-x_{t_r}\cdots x_{\gt_{\ga_r+2}-1}}
	{1-x_{t_r}\cdots x_{\gt_{\ga_r+1}-1}}, \cdots,
\frac{1-x_{t_r}\cdots x_{t_{r+1}-1}}
	{1-x_{t_r}\cdots x_{\gt_{\ga_{r+1}-1}-1}}\right).
\end{multline*}
By Theorem 4.4 and Proposition 5.5 of 
\cite{Zana} $E_{\bfi,\bfj}(\bfx)$ has monodromy
along $\caD_{j0}$ if and only if $t_r\le j\le t_{r+1}-2$ for some $r\ge 
1$.
According to the computation in Proposition 5.5 we
further have that
$$(\gTH(q_{j0})-\id)E_{\bfi,\bfj}(\bfx)
=-2\pi i\sum_{s=j+1}^{t_{r+1}-1}E_{\bfi,\bfj+\bfu_s}(\bfx)$$
which involves only the entries on the $\bfi$-th row. Hence
$$(\gTH(q_{j0})-\id)C_\bfj(\bfx)
=-\sum_{s=j}^{t_{r+1}-1}C_{\bfj+\bfu_s}(\bfx).$$

By similar argument using Proposition 5.4 and 5.5 of \cite{Zana} 
we see that if $t_r=i\le j\le t_{r+1}-2$, $r\ge 1$, then
$$(\gTH(q_{ij})-\id)E_{\bfi,\bfj}(\bfx)
=2\pi iE_{\bfi,\bfj+\bfu_{j+1}}(\bfx)$$
and therefore
$$(\gTH(q_{ij})-\id)C_\bfj(\bfx)
=-C_{\bfj+\bfu_{j+1}}(\bfx).$$

Similarly, thanks to Theorem 5.3 and Proposition 5.5
of \cite{Zana} if $t_r+1\le i\le j=t_{r+1}-1,\ r\ge 0$, then
$$(\gTH(q_{ij})-\id)E_{\bfi,\bfj}(\bfx)
=-2\pi iE_{\bfi,\bfj+\bfu_i}(\bfx).$$
Hence 
$$(\gTH(q_{ij})-\id)C_\bfj=-C_{\bfj+\bfu_i}(\bfx).$$
This completes the proof of the proposition.
\end{proof}

\begin{cor} The monodromy representation of $\caM_{[n]}(\bfx)$
$$\rho_\bfx: \pi_1(S_n,\bfx)\lra \GL_{2^n}(\Z)$$
is unipotent.
\end{cor}
\begin{proof}
Clear.
\end{proof}

\subsection{MHS of multiple logarithms}
Define a meromorphic connection $\nabla$ on the trivial bundle
\begin{equation}\label{trivbun}
\CP^n\times \C^{2^n}\lra \CP^n
\end{equation}
by 
$$\nabla f=df-{\boldsymbol \om} f$$
where $f: S_n\ra \C^{2^n}$ is a section. This connection has regular 
singularities along $\caD_n$ because $\boldsymbol \om$ is 
integrable by \eqref{integrabl} and all the 1-forms in $\boldsymbol \om$
are logarithmic in any compactification of $S_n$. By the explicit
construction of $\boldsymbol \om$ we see immediately that the conditions
(iii) and (iv) of Theorem~\ref{del} are satisfied. Proposition~\ref
{griff} 
further implies that the columns $(2\pi i)^{|\bfj|}C_\bfj(\bfx)$  of 
$\caM_{[n]}(\bfx)$ satisfy $\nabla f=0$ and are therefore flat sections
of \eqref{trivbun}. Even though they are multi-valued, their $\Z$-
linear 
span is well defined thanks to Theorem~\ref{thm:mono}.
Hence $V_{[n]}(\bfx)$ forms a local system over $S_n$.

\begin{defn}
The local system $V_{[n]}(\bfx)$ is called the {\em $n$-tuple logarithm
local system}. 
\end{defn}

To define the MHS on $V_{[n]}$ 
we can define the weight filtration by putting $W_{2k+1}=W_{2k}$ and
$$W_{-2k}V_{[n]}(\bfx)=\langle  (2\pi i)^{|\bfi|}C_\bfi:\ |\bfi|\ge
k\rangle_\Q$$
which is the $\Q$ vector space with basis 
$\{(2\pi i)^{|\bfi|}C_\bfi:\ |\bfi|\ge k\}$.
In particular, $W_{-2k}V_{[n]}(\bfx)=0$ if $k>n$ and
$W_{-2k}V_{[n]}(\bfx)=V_{[n]}(\bfx)$ if $k\le 0$. By regarding 
$e_\bfi$'s as
column vectors one can define the Hodge filtration on 
$V_{[n]}(\bfx)\ot \C=V_{[n],\C}$
as follows:
$$\caF^{-k}V_{[n],\C} :=\langle  e_\bfi:\ |\bfi|\le k\rangle_\C.$$
So in particular, $\caF^{-k}V_{[n],\C}=0$ for $k< 0$ and
$\caF^{-k}V_{[n],\C}=V_{[n],\C}$ for $k\ge n$.

By induction on $n$ and using Lemma~\ref{kthdiag}
it is easy to show that
$$\caF^{-p}\cap W_{-2k}V_{[n],\C}=\begin{cases}
0  &\text{if }p\le k-1 \\
\langle  (2\pi i)^{|\bfi|} e_\bfi:\ k\le |\bfi|\le p\rangle \quad&\text
{if }
k\le p\le n\\
\langle  (2\pi i)^{|\bfi|} e_\bfi:\ k\le |\bfi|\le n\rangle &\text{if } 
p\ge n
\end{cases}.$$
This implies that
$$
\caF^{-p}\gr_{-2k}^W V_{[n],\C}=\begin{cases}
0  &\text{if } p\le k-1\\
W_{-2k}V_{[n],\C}/W_{-2k-1}V_{[n],\C} \qquad&\text{if }p\ge k.
\end{cases}$$
In other words, $\caF^q \gr_{-2k}^W V_{[n],\C}=0$ for $q\ge -k+1$
and $\caF^q\gr_{-2k}^W V_{[n],\C}=\gr_{-2k}^W V_{[n],\C}$
for $q\le -k$. This
means that the Hodge filtration induces a pure HS of
weight $-2k$ on each weight graded piece. Furthermore, it is not hard to
see by checking the powers of $2\pi i$ appearing on the diagonal
of $\caM_{[n]}(\bfx)$ that this induced structure on $\gr_{-2k}^WV_
{[n],\C}$ is
isomorphic to the direct sum of ${n\choose k}$ copies of the Tate 
structure $\Z(k)$ by Lemma~\ref{kthdiag}.

\section{Limit MHS of multiple logarithms}
Let the monodromy of $\caM_{[n]}(\bfx)$ at any subvariety
$\caD$ of $\CP^n$ be given by the matrix $T_\caD$ and the local
monodromy logarithm by $N_\caD=\log T_\caD/2\pi i$. Note
that $T_\caD$ is unipotent so $N_\caD$ is well-defined.

Now let us recall the construction of the unipotent variations
of limit MHS at the ``infinity'' with {\em normal crossing}. 
Let $S$ be a complex manifold of dimension $d$. 
Suppose that $S$ is embedded in $\tilde{S}$,
via the mapping $j$, such that $D=\tilde{S}-S$ is a divisor
with normal crossings. Let $\V$ be any local system  of complex
vector spaces on $S$, and $\caV$ the corresponding vector bundle.
According to  Theorem~\ref{del} by Deligne
there is a canonical extension $\tilde{\caV}$ of $\caV$
over $\tilde{S}$. Moreover, when the local monodromy
is nilpotent $\tilde{\caV}$ is a subsheaf of $j_{\ast} \caV$. The
local picture of $S\subset \tilde{S}$ is 
$(\gD^*)^r\times \gD^{d-r}\subset \gD^d$
where $\gD$ is the unit disk and $\gD^*$ is the punctured one.
We let $t_1,\dots,t_r$ denote the variables on $(\gD^*)^r$, and 
$N_1,\dots, N_r$ the (commuting) local nilpotent 
logarithms of the associated monodromy transformations of the fibre. 
For $z_1,\dots, z_r$ in the upper half-plane, the universal 
covering mapping for $(\gD^*)^r$ is given by 
$$t_j=\exp(2\pi iz_j),\quad j=1,\cdots,r.$$
Let $v_1,\dots,v_m$ be a basis of the multi-valued sections of $\V$
over $(\gD^*)^r\times \gD^{d-r}$, the formula
$$[\tilde{v}_1,\dots,\tilde{v}_m]
=[v_1,\dots,v_m]\exp\left(-\sum_{j=1}^r 2\pi iz_j N_j\right)
=[v_1,\dots,v_m]\prod_{j=1}^r t_j^{-N_j} $$
determines a basis of the sections of $\caV$ over $\gD^d$ and
these provide, by definition, the generators of $\tilde{\caV}$ over 
$\gD^d$.
 
In our situation, although the divisor $D_n$ is not normal crossing
Theorem~\ref{del} is still valid.
We further notice that the image of the global 
holomorphic logarithmic forms in the complex of smooth forms 
on $S$ is independent of the normal crossings compactification 
(see \cite[Prop.~(3.2)]{HM}).
In fact, the forms we are considering lie in the
subcomplex generated by 1-forms of the type $df/f$ where $f$ is a 
rational function. Such forms are automatically logarithmic in
any compactification and therefore our connection is 
automatically regular. Hence the admissibility and the existence
of the limit MHS is an automatic consequence of the admissibility
of our variations restricted to every curve in $S_n$. Moreover,
the pullback of our trivial bundle \eqref{trivbun} restricted
to $S_n$ to $\tilde{S}_n$
is exactly Deligne's canonical extension of \eqref{trivbun},
and the pullbacks of the subbundles $\caF^\bul$ and $W_\bul$ are the
correct extended Hodge and weight subbundles.
Therefore we have

\begin{thm} The $n$-tuple logarithm underlies a good unipotent 
graded-polarizable variation of 
mixed Hodge-Tate structures $(V_{[n]},W_\bul,\caF^\bul)$ over 
$$S_n=\C^n\setminus \Bigl\{\prod_{1\le j \le n}x_j(1-x_j)
\prod_{1\le i<j\le n} \big(1-x_i\dots x_j\big)=0\Bigr\}.$$
with the weight-graded quotients $\gr_{-2k}^W$ 
being given by ${n\choose k}$ copies of the Tate structure $\Z(k)$.
\end{thm}
\begin{proof} It is clearly that all the odd graded weight quotients 
are zero so that we can let the polarizations on the weight graded 
quotients  $\gr_{-2k}^W$ be the ones that give each vector 
$2\pi i e_\bfj$ ($|\bfj|=k$) length 1. 
Then everything is clear except the Griffiths transversality condition. 
But this condition is also satisfied because 
$dC_\bfj={\boldsymbol \om}C_\bfj$ for every $\bfj\in \SSS_n$
by Proposition~\ref{griff}.
\end{proof}

If we want to determine the limit MHS of multiple logarithms
explicitly we can still apply the techniques used
in the normal crossing case. We will carry this out only
for the depth two and three cases. The general picture is
similar but much more complicated.

\subsection{Limit MHS of double logarithm}
First we look at the double logarithm variation of MHS. We have
$$\caM_{1,1}(x,y)=\left[\begin{matrix}
1&\ &\ &\ \\
\gemL_1(y)  & 2\pi i&\ &\ \\
\gemL_1(xy) & 0&2\pi i&\ \\
\gemL_2(x,y)& 2\pi i\gemL_1(x)&
	2\pi iH(x,y)&(2\pi i)^2
\end{matrix}\right]
$$
where $H(x,y)=\gemL_1(y)-\gemL_1(x)-\log x$.

\noindent(i) Let us first try to extend the MHS to the divisor 
$\caD_{10}=\{x=0\}$ along the tangent vector
$\partial /\partial x$. We have
$$T_{\{x=0\}}=\left[\begin{matrix}
1&\ &\ &\ \\
0&1&\ &\ \\
0&0&1&\ \\
0&0&-1&1\end{matrix}\right],\quad
N_{\{x=0\}}=\frac{\log T_{\{x=0\}}}{2\pi i}=\left[\begin{matrix}
0&\ &\ &\ \\
0&0&\ &\ \\
0&0&0&\ \\
0&0&-\frac{1}{2\pi i}&0\end{matrix}\right].$$
Let $\caM_{1,1}(x,y)=[C_0(x,y)\ \cdots C_3(x,y)]$. Define
$$\aligned
\ [s_0\ s_1\ s_2\ s_3]
=&\lim_{t\to 0}\caM_{1,1}(t,y)\left[\begin{matrix}
1&\ &\ &\ \\
0&1&\ &\ \\
0&0&1&\ \\
0&0&\log t/(2\pi i) &1\end{matrix}\right] \\
=&\left[\begin{matrix}
1&\ &\ &\ \\
\gemL_1(y)  & 2\pi i&\ &\ \\
0& 0&2\pi i&\ \\
0& 0& 2\pi i \gemL_1(y) &(2\pi i)^2
\end{matrix}\right].
\endaligned$$
Let $V_{\Q,\{x=0\}}$ be the $\Q$-linear span of $s_0,s_1,s_2,s_3$,
and $V_{\C,\{x=0\}}=\C\ot V_{\Q,\{x=0\}}$. Let $\{e_j:j=0,\cdots, 3\}$ 
be the
standard basis of $\C^4$ where the only nonzero entry of $e_j$ is
at the $(j+1)$st component. 
Then the limit MHS on $\{(x,y):x=0,y\ne 1\}$ along $\partial /\partial x$ 
are given by 
$$((V_{\Q,\{x=0\}}, W_\bul),
(V_{\C,\{x=0\}}, F^\bul))$$
where for $k=0,\dots,3$
\begin{equation}\label{weight}
W_{-2k} V_{\Q,\{x=0\}}=\langle s_k,\dots, s_3\rangle, W_{-2k}=W_{-2k+1}
\end{equation}
and
\begin{equation}\label{hodge}
F^{-k}V_{\C,\{x=0\}}=\langle e_0,\dots, e_k\rangle.
\end{equation}

\noindent(ii) A similar calculation shows that along the tangent vector  
$\partial /\partial x$ the limit MHS on the divisor
$\caD_{11}=\{(1,y): y\ne 1\}$ is the $\Q$-linear span of $s_0,\dots,s_3$
where
$$ [s_0\ s_1\ s_2\ s_3]
=\left[\begin{matrix}
1&\ &\ &\ \\
\gemL_1(y)  & 2\pi i&\ &\ \\
\gemL_1(y)  & 0&2\pi i&\ \\
\gemL_2(1,y)& 0& 2\pi i \gemL_1(y) &(2\pi i)^2
\end{matrix}\right]. $$
It is easy to see by differentiation that
$\gemL_2(1,y)=\big(\gemL_1(y)\big)^2/2.$

\noindent(iii) The extension of MHS to $\caD_{22}=\{(x,1): x\ne 0,1\}$ 
along the tangent vector $\partial /\partial y$ is given by the $\Q$-linear 
span of $s_0,\dots,s_3$ where
$$ 
\ [s_0\ s_1\ s_2\ s_3]
=\left[\begin{matrix}
1&\ &\ &\ \\
0& 2\pi i&\ &\ \\
-Li_1(\frac{x}{x-1})  & 0&2\pi i&\ \\
Li_2(\frac{x}{x-1}) & 2\pi iLi_1(x) & 
	-2\pi i \log\frac{x}{x-1} &(2\pi i)^2
\end{matrix}\right].$$
 
\noindent(iv) Limit MHS on $\caD_{12}=\{(1/y,y): y\ne 0,1\}$ 
along the tangent vector $\partial /\partial x$ is given by the $\Q$-linear 
span of $s_0,\dots,s_3$ where
$$ [s_0\ s_1\ s_2\ s_3]
=\left[\begin{matrix}
1&\ &\ &\ \\
-Li_1(\frac{y}{y-1})& 2\pi i&\ &\ \\
0  & 0&2\pi i&\ \\
-Li_2(\frac{y}{y-1}) & 2\pi i\log\frac{y}{y-1}& 0 &(2\pi i)^2
\end{matrix}\right].$$
  
\noindent(v) $\caD_{10}\cap \caD_{22}=(0,1)$.  From (i)
we see that there are limit MHS on the open set $\caD_{10}\setminus
\{(0,1)\}$ of $\caD_{10}$. We now can easily extend these MHS 
to $(0,1)$ along the vector $\partial/\partial y$ and find  
the limit MHS at $(0,1)$ to be the $\Q$-linear span of $s_0,\cdots, s_3$ 
where
$$[s_0\ s_1\ s_2\ s_3]= \left[\begin{matrix}
1&\ &\ &\ \\
0& 2\pi i&\ &\  \\
0& 0&2\pi i&\ \\
0 & 0& 0 &(2\pi i)^2
\end{matrix}\right].$$
If we start from (iii) and then extend the MHS to $(0,1)$ 
along tangent vector $\partial/\partial x$ we will get the
same limit MHS.

\noindent(vi) $\caD_{11}\cap \caD_{12}=\caD_{12}\cap \caD_{22}
=\caD_{11}\cap \caD_{22}=(1,1)$. We can start from either case (ii)
or (iii) or (iv). Extending the limit MHS of case (ii) we see 
immediately that the along the tangent vector $\partial/\partial y$ 
the limit MHS at $(1,1)$ is given by the $\Q$-linear span of
\begin{equation}\label{basis1}
[s_0\ s_1\ s_2\ s_3]=\left[\begin{matrix}
1&\ &\ &\ \\
0& 2\pi i&\ &\ \\
0& 0&2\pi i&\ \\
E_{4,1} & 0& 0 &(2\pi i)^2
\end{matrix}\right]. 
\end{equation}
If we extend the limit MHS of case (iii) to $(1,1)$ along tangent vector 
$\partial/\partial x$ we find that only the lower left corner entry
is different from the above. Instead of $0$ it is
$$E_{4,1}=\lim_{x\to 1}
Li_2(\frac{x}{x-1})+\frac12 \log^2(1-x)-\log x\log(1-x)
=-Li_2(1)=-\frac{\pi^2}{12},$$
since 
\begin{equation}\label{li21}
Li_2(1-t)+Li_2(1-1/t)+\log^2t/2=0 \forall t\ne 0.
\end{equation}
But if we take $s_0'=s_0-s_3/48$ we get the same basis
as in \eqref{basis1}. The same phenomenon occurs if 
we start from case (iv) and then use tangent vector 
$\partial/\partial y$. 
 
If we extend the limit MHS of (iv) to the point $(1,1)$ along the 
tangent vector $\partial/\partial y$ then we find that
$$E_{4,1}=\lim_{y\to 1}
-Li_2(\frac{y}{y-1})-\frac12 \log^2(1-y)=Li_2(1)=\frac{\pi^2}{12}$$
by taking $t=1-y$ in \eqref{li21}.
Now if we let $s_0'=s_0+\frac{1}{48}s_3$ then we get the same basis
as in \eqref{basis1}.  This phenomenon happens in higher logarithm
cases too.

\subsection{Limit MHS of triple logarithm}
The triple logarithm function $\gemL_3(x,y,z)$ is defined by
(\cite[Example 5.2]{Zana})
$$\aligned
Li_{1,1,1}(x,y,z)=&\int_{(0,0,0)}^{(x,y,z)}
\frac{dz}{1-z}\frac{dy}{1-y} \frac{dx}{1-x}
+\frac{d(yz)}{1-yz}
\left(\frac{dz}{1-z}+\frac{dy}{y(y-1)}\right)
\frac{dx}{1-x}\\
+&\frac{d(yz)}{1-yz}\frac{dx}{1-x}
\left(\frac{dz}{1-z}+\frac{dy}{y(y-1)}\right)
+\frac{dz}{1-z}\frac{d(xy)}{1-xy}
\left(\frac{dy}{1-y}+\frac{dx}{x(x-1)}\right)\\
+&\frac{d(xyz)}{1-xyz}
\left(\frac{dz}{1-z}+\frac{d(xy)}{xy(xy-1)}\right)
\left(\frac{dy}{1-y}+\frac{dx}{x(x-1)}\right)\\
+& \frac{d(xyz)}{1-xyz}
\left(\frac{d(yz)}{1-yz}+\frac{dx}{x(x-1)}\right)
\left(\frac{dz}{1-z}+\frac{dy}{y(y-1)}\right).
\endaligned$$
Set 
$$\tau_{[3]}(2\pi i)=\text{diag}[1,2\pi i,2\pi i,2\pi i,
	(2\pi i)^2,(2\pi i)^2,(2\pi i)^2,(2\pi i)^3],$$
and define the matrix $\caM_{[3]}(x,y,z)\tau_{[3]}(2\pi i)^{-1}$ by
$$ \begin{pmatrix}
1&\ &\ &\ &\ &\ &\ &\  \\ 
\gemL_1(z)& 1 &\ &\ &\ &\ &\ &\  \\ 
\gemL_1(yz) &0&1&\ &\ &\ &\ &\ \\ 
\gemL_1(xyz)& 0 &0&1&\ &\ &\ &\  \\
\gemL_2(y,z)& \gemL_1(y)& H(y,z)&0&1&\ &\ &\  \\ 
\gemL_2(xy,z)& \gemL_1(xy)&0& H(xy,z)&0&1 &\ &\  \\
\gemL_2(x,yz)&0& \gemL_1(x)& H(x,yz)&0&0&1 &\  \\ 
\gemL_3(x,y,z)&\gemL_2(x,y)&H(y,z)\gemL_1(x)& E_{8,4}&
\gemL_1(x)&H(x,y)&H(y,z) &1 
\end{pmatrix},$$ 
where 
$$ E_{8,4}=\gemL_2\Bigl(\frac{1-xy}{1-x},\frac{1-xyz}{1-xy}\Bigr)
	=\gemL_2(y^{-1},x^{-1})-\gemL_2(y^{-1},yz)+
	\log\frac{x-1}{x(1-yz)} \gemL_1(z).$$

\noindent(i) Extension to an open set of $\caD_{10}=\{x=0\}$ along the vector
$\partial/\partial x$. By Theorem~\ref{thm:mono} or direct
computation $T_{\{x=0\}}=I_8-e_{64}-e_{74}-e_{86}$
so
$$t^{-N_{\{x=0\}}}= I_8+\frac{\log t}{2\pi i}(e_{64}+e_{74}+e_{86})
+\log^2t/2(2\pi i)^2e_{84}.$$   
To determine the limit MHS along $\{x=0\}$ we need to find
$g(y,z)=\lim_{t\to 0} I(t)$ where
$$I(t)=\gemL_2\Bigl(\frac{1-ty}{1-t},\frac{1-tyz}{1-ty}\Bigr)
-\log t\log \big(ty(z-1)\big)+\frac{\log^2t}{2}$$
because $\lim_{t\to 0} \log t\log(1-t)=0.$
We see that $I'(t)=f'(t)$ where
$$f(t)=Li_2(1-t)-\log(1-t)\log\frac{y(1-z)}{y-1}
	-Li_2\Bigl(\frac{(1-t)y}{y-1}\Bigr)+Li_2(ty).$$
Thus
\begin{multline*}
g(y,z)=I(1/yz)+f(0)-f(1/yz)
=Li_2(1-yz)+Li_2(1)+\log\frac{y(1-z)}{y-1}\log(1-yz)\\
-Li_2\Bigl(\frac{y}{y-1}\Bigr)+Li_2\Bigl(\frac{1-yz}{z(1-y)}\Bigr)
-Li_2\Bigl(\frac{1}{z}\Bigr)+\log(1-y)\log(yz).
\end{multline*}
It is easy to see by differentiation with respect to $y$ and $g(0,z)
=2Li_2(1)$ that
$$g(y,z)=Li_{1,1}(y,z)-Li_2(1-y)+3Li_2(1).$$
Hence the local system $V_{\Q,\{x=0\}}$ of the limit MHS over 
$\{(0,y,z): y(1-y)(1-z)(1-yz)\ne 0\}$ 
is the $\Q$-linear span of $s_0,\cdots, s_7$ where
$[s_0\ \cdots\ s_7]$ is given by
$$\begin{pmatrix}
1&\ &\ &\ &\ &\ &\ &\  \\ 
\gemL_1(z)& 1 &\ &\ &\ &\ &\ &\  \\ 
\gemL_1(yz) &0&1&\ &\ &\ &\ &\ \\ 
0 & 0 &0&1&\ &\ &\ &\  \\
\gemL_2(y,z)& \gemL_1(y)& H(y,z)&0&1&\ &\ &\  \\ 
0&0&0&\gemL_1(z)-\log y&0&1 &\ &\  \\
0&0& 0&\gemL_1(yz)&0&0&1 &\   \\ 
0&0&0&g(y,z)& 0&\gemL_1(y)& H(y,z)&1 
\end{pmatrix} \tau_{[3]}(2\pi i).$$ 

\noindent(ii) On $\caD_{20}=\{y=0\}$.  
Similar computation as above shows that
the local system $V_{\Q,\{y=0\}}$ of the limit MHS over 
$\{(x,0,z): x(1-x)(1-z)\ne 0\}$ along the vector
$\partial/\partial y$  is the $\Q$-linear span of $s_0,\cdots, s_7$ where
$[s_0\ \cdots\ s_7]$ is given by
$$ \begin{pmatrix}
1&\ &\ &\ &\ &\ &\ &\  \\ 
\gemL_1(z)& 1 &\ &\ &\ &\ &\ &\ \\ 
0& 0&1&\ &\ &\ &\ &\ \\ 
0& 0&0&1&\ &\ &\ &\  \\
0& 0&\gemL_1(z)&0&1&\ &\ &\ \\ 
0& 0&0&\gemL_1(z)-\log x&0&1 &\ &\  \\
0& 0& \gemL_1(x)&-\gemL_1(x)-\log x &0&0&1 &\  \\ 
0& 0&\gemL_1(z)\gemL_1(x)&g(x,z)&\gemL_1(x)&-\gemL_1(x)-\log x &\gemL_1
(z)&1 
\end{pmatrix} \tau_{[3]}(2\pi i),
$$ 
where
$$g(x,z)=Li_2(1)-Li_1(z)(Li_1(x)+\log x)-Li_2(1-x^{-1}).$$

\noindent(iii) On $\caD_{11}=\{x=1\}$. Then 
the local system $V_{\Q,\{x=1\}}$ of the limit MHS over 
$\{(1,y,z): y(1-y)(1-z)(1-yz)\ne 0\}$ along the vector
$\partial/\partial x$ is the $\Q$-linear span of $s_0,\cdots, s_7$ where
$[s_0\ \cdots\ s_7]$ is given by
$$ \begin{pmatrix}
1&\ &\ &\ &\ &\ &\ &\  \\ 
\gemL_1(z)& 1 &\ &\ &\ &\ &\ &\ \\ 
\gemL_1(yz) &0&1&\ &\ &\ &\ &\ \\ 
\gemL_1(yz)& 0 &0&1&\ &\ &\ &\ \\
\gemL_2(y,z)& \gemL_1(y)& H(y,z)&0&1&\ &\ &\ \\ 
\gemL_2(y,z)& \gemL_1(y)&0&H(y,z)&0&1 &\ &\  \\
\gemL_2(1,yz)&0&0&\gemL_1(yz)&0&0&1 &\ \\ 
\gemL_3(1,y,z)&\gemL_2(1,y)&0& g(y,z)&0&\gemL_1(y)&H(y,z)&1 
\end{pmatrix} \tau_{[3]}(2\pi i),
$$ 
where
$$g(y,z)=\gemL_2(y,z)+Li_2\big(1/(1-y)\big).$$
 
\noindent(iv) On $\caD_{22}=\{y=1\}$. The local system 
$V_{\Q,\{y=1\}}$ of 
the limit MHS over $\{(x,1,z): x(1-x)(1-z)(1-xz)\ne 0\}$ 
along the vector $\partial/\partial y$ 
is the $\Q$-linear span of $s_0,\cdots, s_7$ where
$[s_0\ \cdots\ s_7]$ is given by
$$ \begin{pmatrix}
1&\ &\ &\ &\ &\ &\ &\  \\ 
\gemL_1(z)& 1 &\ &\ &\ &\ &\ &\ \\ 
\gemL_1(z) &0&1&\ &\ &\ &\ &\ \\ 
\gemL_1(xz)&0&0&1&\ &\ &\ &\  \\
\gemL_2(1,z)&0&\gemL_1(z)&0&1&\ &\ &\   \\ 
\gemL_2(x,z)& \gemL_1(x)&0&
	H(x,z)&0&1 &\ &\ \\
\gemL_2(x,z)&0& \gemL_1(x)&H(x,z)& 0&0&1 &\ \\ 
\gemL_3(x,1,z)&\gemL_2(\frac{x}{x-1}) &\gemL_1(x)\gemL_1(z) 
	& \gemL_2(1, \frac{1-xz}{1-x})&\gemL_1(x)&
	-\gemL_1(x)-\log x&\gemL_1(z) &1 
\end{pmatrix} \tau_{[3]}(2\pi i).
$$ 

\noindent(v) On $\caD_{33}=\{z=1\}$. This case is the most interesting
because the variation of MHS for $Li_{2,1}$ appears implicitly.
 
The local system $V_{\Q,\{z=1\}}$ 
of the limit MHS over $\{(x,y,1): xy(1-x)(1-y)(1-xy)\ne 0\}$ 
along the vector $\partial/\partial z$ 
is the $\Q$-linear span of $s_0,\cdots, s_7$ where
$[s_0\ \cdots\ s_7]$ is given by
$$ \begin{pmatrix}
1&\ &\ &\ &\ &\ &\ &\  \\ 
0&1 &\ &\ &\ &\ &\ &\ \\ 
\gemL_1(y) &0&1&\ &\ &\ &\ &\ \\ 
\gemL_1(xy)& 0 &0&1&\ &\ &\ &\  \\
Li_2(\frac{y}{y-1}) & \gemL_1(y)& \log\frac{y-1}{y}&0&1&\ &\ &\  \\ 
Li_2(\frac{xy}{xy-1}) &\gemL_1(xy)&0&
	\log \frac{xy-1}{xy}& 0&1 &\ &\  \\
\gemL_2(x,y)&0& \gemL_1(x)& H(x,y)&0&0&1 &\  \\ 
g(x,y)&\gemL_2(x,y)&\log \frac{y-1}{y}\gemL_1(x)& h(x,y)&
\gemL_1(x)&H(x,y)&\log \frac{y-1}{y} &1 
\end{pmatrix} \tau_{[3]}(2\pi i),
$$ 
where
$$g(x,y)=Li_{1,2}\Bigl(\frac{x(y-1)}{xy-1},\frac{y}{y-1}\Bigr) 
+\log(1-xy)Li_2\Bigl(\frac{y}{y-1}\Bigr)$$
and
$$h(x,y)=Li_2\Bigl(\frac{1-xy}{x(1-y)} \Bigr) 
	+H(x,y)\log\frac{xy-1}{xy}.$$
 
We observe that this is {\em essentially} the variation matrix 
$\caM_{1,2}\Bigl(\frac{x(y-1)}{xy-1},\frac{y}{y-1}\Bigr)$.

We omit the following similar cases: 

(vi) On $\caD_{12}=\{xy=1\}$. Extend along the vector
$\partial/\partial x$ or $\partial/\partial y$,   

(vii) On $\caD_{23}=\{yz=1\}$. Extend along the vector
$\partial/\partial y$ or $\partial/\partial z$, 

(viii) On  $\caD_{13}=\{xyz=1\}$. Extend along the vector
$\partial/\partial x$, or $\partial/\partial y$, or $\partial/\partial z$.

\noindent (ix). $\caD_{10}\cap \caD_{20}.$ 
We may start from either case (i) or case (ii). 
Straightforward calculation starting from case (i) shows that 
the extension of the MHS on $\caD_{10}$ to
$\caD_{10}\cap \caD_{20}$ along the vector $\partial/\partial y$ 
is the $\Q$-linear span of $[s_0\ \cdots\ s_7]$ given by
$$ \begin{pmatrix}
1&\ &\ &\ &\ &\ &\ &\  \\ 
\gemL_1(z)& 1 &\ &\ &\ &\ &\ &\ \\ 
0& 0&1&\ &\ &\ &\ &\ \\ 
0& 0&0&1&\ &\ &\ &\  \\
0& 0&\gemL_1(z)&0&1&\ &\ &\  \\ 
0& 0&0&\gemL_1(z)&0&1 &\ &\   \\
0& 0& 0&0 &0&0&1 &\  \\ 
0& 0& 0& 2Li_2(1) & 0&0& \gemL_1(z)&1 
\end{pmatrix} \tau_{[3]}(2\pi i).$$ 
If we start from case (ii) and take the vector $\partial/\partial x$ 
then we will get the same result.

\noindent(x) $\caD_{11}\cap \caD_{22}.$ We may start from
either case (iii) or case (iv). Straightforward calculation 
starting from case (iii) shows that 
along the vector $\partial/\partial y$ 
the limit MHS on $\caD_{11}\cap \caD_{22}$ is the $\Q$-linear span
of $[s_0\ \cdots\ s_7]$ given by
$$ \begin{pmatrix}
1&0&0&0&0&0&0&0 \\ 
\gemL_1(z)& 1&\ &\ &\ &\ &\ &\ &\  \\ 
\gemL_1(z) &0&1&\ &\ &\ &\ &\ &\ \\ 
\gemL_1(z)& 0 &0&1&\ &\ &\ &\ &\ \\
\gemL_2(1,z)&0&\gemL_1(z) &0&1&\ &\ &\  \\ 
\gemL_2(1,z)&0&0&\gemL_1(z) &0&1&\ &\  \\
\gemL_2(1,z)&0&0&\gemL_1(z)&0&0&1 &\  \\ 
\gemL_3(1,1,z)&0&0 &E_{8,4}&0&0 & \gemL_1(z)&1 
\end{pmatrix} \tau_{[3]}(2\pi i),
$$ 
where $E_{8,4}= \gemL_2(z)+2Li_2(1).$
If we start from case (iv) then we find that $E_{8,4}= \gemL_2(z)$
and therefore we get the same limit MHS on
$\caD_{11}\cap \caD_{22}$ along vector $\partial/\partial x$.

By similar computation we can determine the limit MHS on 
the intersections of any two of the irreducible components $\caD_{ij}$
along any vector.
Finally, at all the of the following four points: 
$(0,0,1)$, $(1,0,1)$, $(0,1,1)$ and $(1,1,1)$ 
we find without much difficulty that 
the columns of the matrix $\tau_{[3]}(2\pi i)$ provide us
$s_0,\dots, s_7$ for the limit MHS along vectors $\partial/\partial x$,
or $\partial/\partial y$, or $\partial/\partial z$.

 From all the above examples we want to make the following
\begin{conj}
The variations of mixed Hodge-Tate structures related to
any multiple polylogarithm can be produced as the variations
of some limit mixed Hodge-Tate structures related to some suitable
choice of multiple logarithm.
\end{conj}

\section{Double polylogarithm variations of MHS}\label{wt3var}
One can similarly generalize the above theory to multiple
polylogarithms. One knows that
on $\C^\times \setminus\{1\}$ the matrix $\caM_n(x)$
\begin{equation*}
\left[\begin{matrix} 1 \\
Li_1(x) & 1 & \ \\
Li_2(x) & \log x & 1 & \\
Li_3(x) & \frac {\log^2 x}2  & \log x  & 1 \\
\vdots  & \vdots            & \vdots   &  \ddots &  \ddots  \\
Li_{n-1}(x) & \frac {\log^{n-2}x}{(n-2)!} & \frac {\log^{n-3}x}{(n-3)!} 
 &\cdots & \log x & 1 \\
Li_n(x) & \frac {\log^{n-1}x}{(n-1)!} & \frac {\log^{n-2}x}{(n-2)!} 
& \cdots& \frac {\log^2 x}2  & \log x  & 1
\end{matrix}\right]\text{diag}\Bigl[1,2\pi i,\dots,(2\pi i)^n\Bigr]
\end{equation*}
provides a variation of mixed Hodge-Tate structures related
to the classical $n$-logarithm (cf. \cite{Hain}). To be more
precise, in the definition of $Li_m(x)$ and $\log^m(x)/m!$ above
we actually fixed a path $p$ from $0$ to $x$ and a path $q$ from
1 to $x$ (both independent of $m$) and set
$$
Li_m(x)=\int_p \frac{dt}{1-t}\ub{ \frac{dt}{t}\cdots\frac{dt}{t}}_{
m-1\text{ times}},\qquad
\frac{\log^m(x)}{m!}=\int_q\ub{ \frac{dt}{t}\cdots\frac{dt}{t}}_{
m\text{ times}}.
$$

\subsection{Double polylogarithms of weight 3}
There are only four multiple polylogarithms of weight 3.
Having dealt with $Li_3$ and $Li_{1,1,1}$ we now turn to $Li_{1,2}$
and $Li_{2,1}$.
\begin{thm}\label{wt3thm}
Each of the weight three depth two multiple polylogarithms
underlies a good variation of mixed Hodge-Tate structures over
$S_2=\C^2\setminus\{xy(1-x)(1-y)(1-xy)=0\}$.
For $Li_{2,1}$ the graded weight quotients are $\Z(0)$, $\Z(1)\oplus\Z
(1)$,
$\Z(2)\oplus\Z(2)$, and $\Z(3)$. For $Li_{1,2}$ they
are $\Z(0)$, $\Z(1)\oplus\Z(1)$, $\Z(2)\oplus\Z(2)\oplus\Z(2)$, and $\Z
(3)$.
\end{thm}

\begin{proof} Let
$\tau_{2,1}(\gl)=\text{diag}\big[1,\gl,\gl,\gl^2,\gl^2,\gl^3\big]$.
We define the multi-valued matrix function over $S_2$
$$\caM_{2,1}(x,y)=\left[\begin{matrix}
1&0&0&0&0&0\\
Li_1(y)  & 1&0&0&0&0\\
Li_1(xy) & 0&1&0&0&0\\
Li_{1,1}(x,y)& Li_1(x) &Li_1\big(\frac{1-xy}{1-x}\big)&1&0&0\\
Li_2(xy) & 0&\log (xy)&0&1&0\\
Li_{2,1}(x,y) & Li_2(x)& f(x,y)&\log x & Li_1(y)  &1\\
\end{matrix}\right]\tau_{2,1}(2\pi i)
$$
where
$$f(x,y)=-\int_{a_1}^1 \frac{dt}{t}\frac{dt}{t-a_2} =
Li_2(x^{-1})-Li_2(y)+\log(xy)Li_1(y).$$
The columns of $\caM_{2,1}(x,y)$ form the fundamental solutions of the
differential equation over $S_2$
\begin{equation*}
d \gl =
\left[\begin{matrix}
0 &0 &0 &0 &0 &0 \\
dLi_1(y) &0 &0 &0 &0 &0 \\
dLi_1(xy)&0 &0 &0 &0 &0 \\
0 & dLi_1(x)  & dLi_1\big(\frac{1-xy}{1-x}\big)&0 &0 & 0 \\
0 & 0  & d\log(xy) & 0 &0 & 0 \\
0 & 0  & 0  & d\log x & dLi_1(y)  & 0 \\
\end{matrix}\right] \gl
\end{equation*}
Let $1\le i\le j\le 2$ and $q_{ij}\in \pi_1(S_2,\bfx)$ 
(resp. $1\le j\le 2$ and $q_{j0}$) be
a loop in $S_2$ turning around the irreducible component $\caD_{ij}$
counterclockwise only once such that
$\int_{q_{ij}} d\log(1-x_i\dots x_j)=-2\pi \sqrt{-1}$ 
(resp. $\int_{q_{j0}}d\log x_j= 2\pi i$). 
Let $e_{st}$ be the matrix with  1 at $(s,t)$-th
entry and 0 elsewhere. Observe that if $q_{i\infty}$ is a a loop in 
$S_2$  turning around $x_i=\infty$ only once then
$q_{i\infty}=-q_{i0}+q_{ii}$. By simple computation we see that
the monodromy representation $\rho:\pi_1(S_2,\bfx)\to \GL_6(\Q)$
is given as follows:
$$\aligned
M(q_{10})=&I-e_{43}+e_{53}+e_{64}\\
M(q_{20})=&I+e_{63}\\
M(q_{11})=&I+e_{42}-e_{43}\\
M(q_{22})=&I+e_{21}+e_{43}+e_{65}\\
M(q_{12})=&I+e_{31}
\endaligned$$

We can now easily define the weight and Hodge filtrations, 
determine the MHS over $S_2$ and compute the limit MHS at the 
``infinity''.
This proves the theorem for $Li_{2,1}$. 

To deal with the multiple polylogarithm $Li_{1,2}(x,y)$ we set
$$\tau_{1,2}(\gl)=\text{diag}\big[1,\gl,\gl,\gl^2,\gl^2,\gl^2,\gl^3\big]
$$
and define the
multi-valued matrix function $\caM_{1,2}(x,y)$ over $S_2$ as
$$ \left[\begin{matrix}
1&0&0&0&0&0&0\\
Li_1(y)  & 1&0&0&0&0&0\\
Li_1(xy) & 0&1&0&0&0&0\\
Li_{1,1}(x,y)& Li_1(x) &Li_1\big(\frac{1-xy}{1-x}\big)&1&0&0&0\\
Li_2(y) &\log (y)&0&0&1&0&0\\
Li_2(xy)& 0&\log (xy)&0&0&1&0\\
Li_{1,2}(x,y) & Li_1(x)\log(y)& g(x,y)&\log y& Li_1(x)& -Li_1(x^{-1})  
&1\\
\end{matrix}\right]\tau_{1,2}(2\pi i)$$
where
$$g(x,y)=-\int_{a_1}^1 \frac{dt}{t-a_2} \frac{dt}{t}
=Li_2(y)-Li_2(x^{-1})-\log(xy)Li_1(x^{-1}).$$
The columns of $\caM_{2,1}(x,y)$ form the fundamental solutions of the
differential equation over $S$
\begin{equation*}
d \gl =
\left[\begin{matrix}
0 &0 &0 &0 &0 &0 \\
dLi_1(y) &0 &0 &0 &0 &0 \\
dLi_1(xy)&0 &0 &0 &0 &0 \\
0 & dLi_1(x)  & dLi_1\big(\frac{1-xy}{1-x}\big)&0 &0 & 0 \\
0 & d\log(y)   & 0 & 0 &0 & 0 \\
0 & 0  & d\log(xy) & 0 &0 & 0 \\
0 & 0  & 0  & d\log(y) & dLi_1(x)  & -dLi_1(x^{-1}) & 0 \\
\end{matrix}\right] \gl
\end{equation*}
The monodromy representation $\rho:\pi_1(S_2,\bfx)\to \GL_7(\Q)$
is given as follows:
$$\aligned
M(q_{10})=&I-e_{43}+e_{63}-e_{76}\\
M(q_{20})=&I+e_{52}+e_{63}+e_{74}\\
M(q_{11})=&I+e_{42}-e_{43}+e_{75}-e_{76}\\
M(q_{22})=&I+e_{21}+e_{43}\\
M(q_{12})=&I+e_{31}
\endaligned$$

We can now determine the MHS over $S_2$ and compute the limit 
MHS at the ``infinity'' as before.
This proves the theorem for $Li_{2,1}$.
\end{proof}

\subsection{Some open problems}
It seems very difficult to write down explicitly the variation matrix
associated with the general multiple polylogarithm
$Li_{m_1,\dots,m_n}(\bfx)$. However, the following  general result 
must be true:
\begin{quote} {\em The multiple polylogarithm
$Li_{m_1,\dots,m_n}(\bfx)$ underlies a good unipotent 
graded-polarizable variation of 
mixed Hodge-Tate structures $(V_{m_1,\dots,m_n},W_\bul,\caF^\bul)$ over 
$$S_n=\C^n\setminus \Bigl\{\prod_{i=1}^n x_i(1-x_i)
\prod_{1\le i<j\le n} \big(1-x_i\dots x_j\big)=0\Bigr\}$$
with the weight-graded quotients $\gr_{-2k}^W$ 
being given by $c_k$ copies of the Tate structure $\Z(k)$ which
are nonzero only for $0\le k\le K$.}
\end{quote}
Here $c_k$ is the number of different ways to pick ordered $(k+2)$-tuples
$(b_{\ga_0},\dots,b_{\ga_{k+1}})$ from the ordered numbers
$(b_0,\dots,b_{K+1})$  in the following tableau where $a_1,\dots,a_n$ are
nonzero
\begin{equation}\label{tabel}\tag{$\ast$}
\Big|\ b_0 \Big|\ \cdots \ \Big|\ b_{K+1}\ \Big|=
\Big|\ 0 \ \Big|\ a_1 \ \Big|\ 
\underbrace{0 \ \Big|\ \cdots \ \Big|\ 0}_{a_1-1 \text{ times}} \ \Big|
\ \ a_2\ \Big|\ \cdots \cdots  \ \Big|\ a_n \ \Big|\ 
\underbrace{0\ \Big|\ \cdots\ \Big|\ 0}_{a_n-1\text{ times}}\ \Big|\ 1\ 
\Big|\ 
\end{equation} 
such that all of the following conditions are satisfied:
\begin{quote} 
(i) $\ga_0=0$,

(ii) $\ga_{k+1}=K+1$,

(iii) For all $0\le i\le k$, either $\ga_{i+1}=\ga_i+1$ or at least 
one of $b_{\ga_i}$ and $b_{\ga_{i+1}}$ is nonzero,
\end{quote}

It is apparent that 
$$
c_k\ge d_k(m_1,\dots,m_n)
=\sum_{\substack{k_1+\dots+k_n=k\\ 0\le k_i\le m_i}} 1.$$
Each term in the sum corresponds to the following choice:
for every $i=1,\dots,n$, choose $k_i$ $0$'s immediately
after $a_i$.

\begin{eg} By the definition, we always have $c_0=c_K=1$.
When $m_1=\cdots=m_n=1$ tableau \eqref{tabel} becomes
$$\Big|\ b_0 \Big|\ \cdots \ \Big|\ b_{n+1}\ \Big|=\Big|\ 0\ \Big|\ 
\underbrace{1\ \Big|\ \cdots\ \Big|\ 1}_{n+1\text{ times}}\ \Big|$$
Because $b_0$ and the last $b_{n+1}$ is always picked, $c_k$ is the 
number of
ways to choose $k$ elements from the set $\{b_1,\dots,b_n\}$, i.e., 
$c_k={n\choose k}$.
\end{eg}

For ease of statement let us put a box $\boxed{\phantom{\cdot}}$ on a 
number whenever we choose it. 
\begin{eg} Let's look at $Li_{1,2}$. We have the following six 
nontrivial ways to
put boxes on $\big|\ 0\ \big|\ a_1\ \big|\ a_2\ \big|\ 0\ \big|\ 1\ 
\big|$:
$$\begin{array}{rlll}
\ & \text{(1)\ }\Big|\boxed{0}\Big|\boxed{\phantom{0}}\hskip-9.4pt a_1
\Big|\
	\phantom{0}\hskip-6pt a_2\Big|\ 0\ \Big|\boxed{1}\Big|\qquad &
\text{(2)\ }\Big|\boxed{0}\Big|\ \phantom{0}\hskip-6pt a_1\Big|
	\boxed{\phantom{0}}\hskip-9.3pt a_2\Big|\ 0\ \Big|\boxed{1}
\Big| \qquad &
\text{(3)\ }\Big|\boxed{0}\Big|\boxed{\phantom{0}}\hskip-9.4pt a_1\Big|
	\boxed{\phantom{0}}\hskip-9.3pt a_2\Big|\ 0\ \Big|\boxed{1}
\Big|  \\
\ \\
\ &\text{(4)\ }\Big|\boxed{0}\Big|\boxed{\phantom{0}}\hskip-9.4pt a_1
\Big|\ 
	\phantom{0}\hskip-6pt a_2\Big|\boxed{0}\Big|\boxed{1}\Big| 
\qquad &
\text{(5)\ }\Big|\boxed{0}\Big|\ \phantom{0}\hskip-6pt a_1\Big|
	\boxed{\phantom{0}}\hskip-9.3pt a_2\Big|\boxed{0}\Big|\boxed{1}
\Big|\qquad &
\text{(6)\ }\Big|\boxed{0}\Big|\boxed{\phantom{0}}\hskip-9.4pt a_1\Big|
	\boxed{\phantom{0}}\hskip-9.3pt a_2\Big|\boxed{0}\Big|\boxed{1}
\Big| 
\end{array}$$
Thus $c_0=c_3=1$, $c_1=2$ and $c_2=3$.

However, for $Li_{2,1}$ we have altogether only six ways to do this:
$$\begin{array}{rlll}
\ & \text{(1)\ }\Big|\boxed{0}\Big|\ \phantom{0}\hskip-6pt a_1\Big|
	\ 0\ \Big|\phantom{0}\hskip-6pt a_2 \Big|\boxed{1}\Big| \qquad &
\text{(2)\ }\Big|\boxed{0}\Big|\boxed{\phantom{0}}\hskip-9.4pt a_1
\Big|\ 0\ \Big|\
	\phantom{0}\hskip-6pt a_2\Big|\boxed{1}\Big|\qquad &
\text{(3)\ }\Big|\boxed{0}\Big|\ \phantom{0}\hskip-6pt a_1\Big|\ 0\ 
\Big|
	\boxed{\phantom{0}}\hskip-9.3pt a_2\Big|\boxed{1}\Big|\\
\ \\
\ &\text{(4)\ }\Big|\boxed{0}\Big|\boxed{\phantom{0}}\hskip-9.4pt a_1
\Big|
	\boxed{0}\Big|\ \phantom{0}\hskip-6pt a_2\Big|\boxed{1}\Big| 
\qquad &
\text{(5)\ }\Big|\boxed{0}\Big|\boxed{\phantom{0}}\hskip-9.4pt a_1
\Big|\ 0\ 
	\Big|\boxed{\phantom{0}}\hskip-9.3pt a_2\Big|\boxed{1}\Big| 
\qquad &
\text{(6)\ }\Big|\boxed{0}\Big|\boxed{\phantom{0}}\hskip-9.4pt a_1
\Big|\boxed{0}\Big|
\boxed{\phantom{0}}\hskip-9.3pt a_2\Big|\boxed{1}\Big| 
\end{array}$$
Thus $c_0=c_3=1$, $c_1=c_2=2$.
\end{eg}

We now can generalize Theorem \ref{wt3thm} to
\begin{thm}
The double polylogarithm $Li_{r,s}$ underlies a good unipotent
graded-polarizable variation
of a mixed  Hodge-Tate structure with the graded weight piece
$\gr_{-2k}^W$ being direct sums of $c_k$ copies of $\Z(k)$ where
$$ c_k=\begin{cases}
d_k(r,s)+1 \quad &\text{if } r\ne k=s,\\
d_k(r,s)&\text{otherwise},
\end{cases}$$
and
$$ 
d_k(r,s)=\begin{cases}
0 &\text{if } k< 0\text{ or } k>r+s ,\\
k+1 \quad &\text{if } 0\le k\le \min\{r,s\}, \\
\min\{r,s\}+1  \quad &\text{if } \min\{r,s\}\le k\le \max\{r,s\}, \\
r+s+1-k  \quad &\text{if }  \max\{r,s\}\le k\le r+s.
\end{cases}
$$
\end{thm}
 
Among all the double polylogarithms the homogeneous one $Li_{r,r}(x,y)$
behaves most regularly. It satisfies $c_0=c_{2r}=1,c_1=c_{2r-1}=2,
\dots, c_{r-1}=c_{r+1}=r, c_r=r+1.$

In general, as we remarked at the beginning of this section,
the multiple polylogarithm $Li_{m_1,\dots,m_n}(\bfx)$ underlies a good
variation of mixed  Hodge-Tate structures with the graded weight piece
$\gr_{-2k}^W$ being direct sums of $c_k$ copies
of $\Z(k)$ for some positive integer $c_k$.
It is clear that $c_k\ge d_k(m_1,\dots,m_n)$ and
$c_k(1,\dots,1,m_n)=d_k(m_1,\dots,m_n)$. 
It would be very interesting to solve the following
\begin{prob}\label{probs} 
\begin{description}
\item{(1)}
 Find a closed formula for $c_k$ depending only on $m_1,\dots,m_n$ and 
$k$.

\item{(2)} 
Determine the variation matrix $\caM_{m_1,\dots,m_n}(\bfx)$ explicitly.

\item{(3)}  Determine the connection matrix $\boldsymbol \om$ 
explicitly.

\item{(4)}  Determine the monodromy actions explicitly.
\end{description}
\end{prob}

\section{Single-valued version of multiple polylogarithms}
If part (2) of Problem \ref{probs} is solved then 
following an idea of Beilinson and Deligne \cite{BD} as given
in \cite{Bl} one can
easily discover the single-valued version of 
$Li_{m_1,\dots,m_n}(x_1,\dots,x_n)$ which we denote by
$\caL_{m_1,\dots,m_n}(x_1,\dots,x_n)$ which should be a
real analytic function. In what follows we outline the procedure for
multiple logarithms only.

\subsection{General procedure for producing single-valued \\
multiple logarithms}
For any $n\ge 2$ let $L_{[n]}=L_{[n]}(\bfx)=[C_\zero \ \dots \ C_\one]$ 
be the matrix with $2^n$ columns $C_\bfj$ ($\bfj\in\caS_n$) as before 
and
$\caM_{[n]}=\caM_{[n]}(\bfx)=L_{[n]}(\bfx)\tau_{[n]}(2\pi i)$ where
$$\tau_{[n]}(\gl)=\text{diag}\big[\gl^{|\bfj|}\big]_{\bfj\in\caS_n}.$$
Define the matrix
$$B_{[n]}=\tau_{[n]}(i)\caM_{[n]} \ol{\caM}_{[n]}^{-1}\tau_{[n]}(i)$$
where $\ol{\caM}_{[n]}$ is the complex conjugation of $\caM_{[n]}$.
 From our calculation of the monodromy we see that $B$ is a
single-valued matrix function defined over $S_n$. Moreover 
$$\ol{B}_{[n]}^{\phantom{\frac12}}= B_{[n]}^{-1}$$
since $\ol{\tau_{[n]}(i)}=\tau_{[n]}(i)^{-1}$.
Now that $B_{[n]}=I+N$  with $I$ the identity matrix and $N$ a nilpotent
matrix we see that $\log B$ is well defined and satisfies
$$\ol{\log B_{[n]}}=-\log B_{[n]},$$
namely, $\log B_{[n]}$ is a pure imaginary matrix. Then we define $-1/
(2i)$
times the lower left corner entry of $\log B$ to be $\caL_{[n]}(\bfx)$
which is a single-valued real analytic version of the multiple
logarithm $\gemL_n(\bfx)$.
\begin{rem} Our method is slightly different from that in \cite{BD}.
In fact when we are in the polylogarithm case
the matrix $B$ constructed as above is the conjugate of the one in 
\cite{BD} by $\tau(i)$.
\end{rem}

\subsection{Single-valued double logarithms}
We have seen that
$$L_{1,1}(x,y)=\left[\begin{matrix}
1  & & & \\
Li_1(y)  & 1 & & \\
Li_1(xy) &  & 1 & \\
Li_{1,1}(x,y) &  Li_1(x) & \log \frac{x-1}{x(1-y)}  &1
\end{matrix}\right]
\text{ and } 
\tau_{1,1}(\gl)=\left[\begin{matrix}
1  & & & \\
  &\gl & & \\
  & & \gl& \\
  & & &   \gl^2\end{matrix}\right].
$$
Let $B_{1,1}(x,y)=\tau_{1,1}(i)L_{1,1}(x,y)\tau_{1,1}(-1)
\ol{L_{1,1}(x,y)}^{-1}\tau_{1,1}(i)$. Then $B_{1,1}(x,y)$ is unipotent
and single-valued. An easy calculation shows
$$\log B_{1,1}(x,y)=\left[\begin{matrix}
0  & & & \\
-2i\log|1-y|  & 0 & & \\
-2i\log|1-xy| &  & 0 & \\
-2i\caL_{1,1}(x,y)
 &  -2i\log|1-x| &  2i\log\bigl|\frac{x-1}{x(1-y)}\bigr| &0
\end{matrix}\right]
$$
where
\begin{equation}\label{llsing}
\caL_{1,1}(x,y)=\im\big(Li_{1,1}(x,y)\big)
-\arg(1-y)\log|1-x|+\arg(1-xy)\log\Bigl|\frac{x-1}{x(1-y)}\Bigr|
\end{equation}
is the single-valued real analytic version of $Li_{1,1}(x,y)$.

By differentiation it is easy to check that
$$Li_{1,1}(x,y)
=Li_2\Bigl(\frac{xy-y}{1-y}\Bigr)-Li_2\Bigl(\frac{y}{y-1}\Bigr)-Li_2
(xy).$$
So by using the single-valued dilogarithm function
$\caL_2(z)=\im\big(Li_2(z)\big)+\arg(1-z)\log|z|$
we can also recover \eqref{llsing} as
\begin{equation}\label{llsing2}
\caL_{1,1}(x,y)=\caL_2\Bigl(\frac{xy-y}{1-y}\Bigr)
-\caL_2\Bigl(\frac{y}{y-1}\Bigr)-\caL_2(xy).
\end{equation}
This function satisfies the functional equations
$$\caL_{1,1}(x,y)
=-\caL_{1,1}\Bigl(1-x,\frac{y}{y-1}\Bigr)$$
by the functional equations $\caL_2(x)=-\caL_2(1-x)=-\caL_2(1/x).$

\subsection{Single-valued double polylogarithms $\caL_{1,2}$ and
$\caL_{2,1}$}
By \cite{Zag} a single-valued version of $Li_3(x)$ can be defined as
\begin{equation}\label{3log}
\caL_3(z)=\re\big(Li_3(z)\big)-\log|z|\re\big(Li_2(z)\big)
-\frac13(\log|z|)^2\log|1-z|.
\end{equation}

We now look at $Li_{2,1}(x,y)$ and $Li_{1,2}(x,y)$. By the procedure 
outlined in the first section of this chapter we find that
the single-valued version of $Li_{1,2}(x,y)$ is 
\begin{multline*}
\caL_{1,2}(x,y)=\re Li_{1,2}(x,y)
-\arg(1-xy)\big[\caL_2(x)+\caL_2(y)\big]+\log|1-x|\re Li_2(y)\\
-\log|y|\re Li_{1,1}(x,y)-\log|1-x^{-1}|\re Li_2(xy)
-\frac13\log|xy^2|\log|1-xy|\log\Bigl|1-x^{-1}\Bigr|\\
+\frac13\log|y|\big(2\log|1-y|\log|1-x|+\log|1-xy|\log|x(1-y)|\big).
\end{multline*}
The single-valued version of $Li_{2,1}(x,y)$ is
\begin{multline*}
\caL_{2,1}(x,y)=\re Li_{2,1}(x,y)
+\arg(1-xy)\big[\caL_2(x)+\caL_2(y)\big]-\arg(1-y)\caL_2(x)\\
+\log|1-y|\re Li_2(xy)-\log|x|\re Li_{1,1}(x,y)
+\frac13\log|1-y|\log|xy|\log|1-xy|\\
+\frac13\log|x|\Bigl[\log|1-y|\log|1-x|
+\log|1-xy|\log\Bigl|\frac{x(1-y)}{1-x}\Bigr|\Bigr].
\end{multline*} 

Using the single-valued versions of dilogarithm $\caL_2(z)$
and trilogarithm $\caL_3(z)$ we can express $\caL_{2,1}(y,x)$
by the trilogarithms
\begin{multline*}
\caL_{2,1}(y,x)=\caL_3(1-xy)+\caL_3(1-x)-\caL_3\Bigl(\frac{1-x}{1-xy}
\Bigr)
-\caL_3(y)+\caL_3\Bigl(\frac{y-xy}{1-xy}\Bigr)-\caL_3(1),
\end{multline*}
where $\caL_3$ is the single-valued trilogarithm given by \eqref{3log}.
This follows from the relation (see \cite{Zadv}) first discovered by 
Zagier
after Goncharov's conviction that such identity should exist:
\begin{multline*}
Li_{2,1}(y,x)=Li_3(1-xy)+Li_3(1-x)-Li_3\Bigl(\frac{1-x}{1-xy}\Bigr)
-Li_3(y)+Li_3\Bigl(\frac{y-xy}{1-xy}\Bigr)-Li_3(1)\\
-\log(1-xy)\big(Li_2(1)+Li_2(1-x)\big)
-\log\Bigl(\frac{1-x}{1-xy}\Bigr)Li_2(y)
+\frac12\log(y)\log^2(1-xy).
\end{multline*}
By straightforward computation we further discover the
following interesting formula:
$$\caL_{1,2}(x,y)+\caL_{2,1}(y,x)+\caL_3(xy)=0.$$
One should compare this with 
$$
Li_{1,2}(x,y)+Li_{2,1}(y,x)+Li_3(xy)=-\log(1-x)Li_2(y).$$
Finally we find the interesting identity
$$\aligned
 \caL_{1,1,1}(x,y,z)=&
\caL_3\Bigl(\frac{(y-1)(1-xyz)}{y(1-x)(1-z)}\Bigr) 
+\caL_3\Bigl(\frac{y}{y-1}\Bigr)
+\caL_3(xy)
-\caL_3\Bigl(\frac{1-xyz}{1-x}\Bigr)\\
\ &-\caL_3\Bigl(\frac{1-xyz}{xy(1-z)}\Bigr)
-\caL_3\Bigl(\frac{y-yz}{y-1}\Bigr) 
-\caL_3\Bigl(\frac{y-xy}{y-1}\Bigr)
+\caL_3(1-x).
\endaligned$$
We remind the readers that such identities in higher weight
cases do not exist in general. For example, $\caL_{2,2}(x,y)$
cannot be expressed by only tetralogarithms $\caL_4$.

\subsection{A problem of multiple Dedekind zeta values}
In general there should exist single-valued real analytic version
of the multiple polylogarithm $Li_{m_1,\dots,m_n}(\bfx)$
which we denote by $\caL_{m_1,\dots,m_n}(\bfx)$. For $m_n\ge 2$
the value of this function when $|x_i|\le 1$ is given by the power
series expansion \eqref{pexp}. We end our paper
by stating a generalized Zagier conjecture about special values of
Dedekind zeta function over number fields. 

Denote by $O_F$ the ring of
integers of a number field $F$ and $I_F$ the set of integral
ideals of $O_F$. Let $N$ be the norm from $F$ to $Q$.
Then we define the multiple Dedekind zeta function of depth $d$ over 
$F$ as
$$\zeta_F(s_1,\dots,s_d)=\sum_{\aligned \gemn_1,\dots,\gemn_d\in& O_F\\
N(\gemn_1)<\cdots < &N(\gemn_d)\endaligned} N(\gemn_1)^{-s_1}\cdots
N(\gemn_d)^{-s_d}.$$
This function is well defined for $\re(s_1)>0,\dots,\re(s_{d-1})>0,\re
(s_d)>1$.

\begin{prob} For any integers $m_1,\dots,m_{d-1}\ge 1$ and $m_d\ge 2$, 
is there an expression of
$\zeta_F(m_1,\dots,m_d)$ in terms of a determinant of $\caL_
{m_1,\dots,m_d}$ evaluated at $F$ rational points
up to some factors determined only by
the number field $F$ (such as the discriminant, the number of real and
complex embeddings, etc.)?
\end{prob} 

When $F=\Q$ the problem has an easy answer: 
$$\zeta_\Q(m_1,\dots,m_d)=\caL_{m_1,\dots,m_d}(1,\dots,1).$$

\begin{rem} H.~Gangl kindly informed the author of
the preprint \cite{Wo} of Wojtkowiak in which conjectures generalizing
Zagier's are also considered. 
\end{rem}

\bigskip

\noindent
{\em Address:} Department of Mathematics, University of Pennsylvania, 
PA 19104, USA

\noindent
{\em Email:} jqz@math.upenn.edu

\end{document}